\documentclass[12pt,english,a4paper]{smfart}

\usepackage[utf8]{inputenc}
\usepackage{graphicx}
\usepackage[T1]{fontenc}
\usepackage[english]{babel}
\usepackage{lmodern}
\usepackage{smfthm}
\usepackage[headings]{fullpage}
\usepackage{amssymb}
\usepackage{amsmath}
\usepackage[all]{xy}
\usepackage[colorlinks=true, allcolors=magenta]{hyperref}
\usepackage{url}
\usepackage{tikz}
\usepackage{natbib}
\theoremstyle{plain}
\newtheorem{thm}{Theorem}
\newtheorem{lemma}[thm]{Lemma}
\newtheorem{cor}[thm]{Corollary}
\newtheorem{pro}[thm]{Proposition}
\theoremstyle{definition}

\newtheorem{dfn}[thm]{Definition}

\def\Gal{{\mathrm{Gal}}}

\title{Multivariable $(\varphi, N)$-modules and $(\varphi, N_\nabla)$-modules}
\author{Csaba Anderlik}

\begin{document}
\begin{abstract}
    The goal of this article is to introduce the category of effective filtered multivariable $(\varphi, N)$-modules and the category of multivariable $(\varphi, N_\nabla)$-modules of finite $E_{\Delta}$-height over the rigid-analytic open polydisk $\mathcal{O}_{\Delta}$, and to establish an exact tensor equivalence between them. This generalizes Kisin's analytic classification of filtered $(\varphi, N)$-modules and provides the rigid-analytic foundation for multivariable $p$-adic Hodge theory.
\end{abstract}
\maketitle

\tableofcontents

\section*{Introduction}
\subsection*{Overview}
The development of $p$-adic Hodge theory has been profoundly shaped by the distinct rational and integral frameworks established in the foundational works of Berger and Breuil. Breuil laid the groundwork for integral $p$-adic Hodge theory through his extensive study of $p$-divisible groups, finite flat group schemes, and strongly divisible modules (\cite{Br4,Br1,Br2,Br3}). In contrast, Berger's articles (\cite{Be3,Be1,Be2}) advanced $p$-adic Hodge theory by translating $p$-adic Galois representations into the language of $(\varphi, \Gamma)$-modules and $p$-adic differential equations.

Building upon this rational analytic framework, Kisin (\cite{Kis1}) established a powerful new classification of crystalline and semi-stable representations. A crucial step in Kisin's method was connecting Fontaine's filtered $(\varphi, N)$-modules (\cite{Fo2}) directly to rigid-analytic differential equations. To achieve this, Kisin introduced the category of $(\varphi, N_\nabla)$-modules over $\mathcal{O}$, the ring of rigid-analytic functions on the open unit disk. By demonstrating that the category of effective, filtered $(\varphi, N)$-modules is equivalent to the category of $(\varphi, N_\nabla)$-modules over $\mathcal{O}$, Kisin effectively bridged linear-algebraic period data with differential equations over a rigid-analytic space. 

To define the category precisely, let $k$ be a perfect field of characteristic $p$, let $W:=W(k)$ the ring of Witt vectors, and let $K_{0}:= W\left[\frac{1}{p}\right]$ be the fraction field of $W$. Let $K/K_{0}$ is a finite, totally ramified field extension of $K_{0}$, and let $\pi\in K$ be a uniformizer with $E(u) \in K_{0}[u]$ Eisenstein polynomial. Let $\overline{K}$ be an algebraic closure of $K$, and define a sequence of elements in $\overline{K}$ for $n \in \mathbb{Z}^{\geq 0}$, $\pi_{0}:=\pi$, $\pi_{n+1}^p=\pi_{n}$. Denote $K_{n+1}:= K(\pi_{n})$. Setting $\mathfrak{S}:=W[[u]]$, we define $\widehat{\mathfrak{S}_{n}}$ as a completion of $K_{n+1}\otimes_{W} \mathfrak{S}$ at the maximal ideal $(u-\pi_{n})$. 

Let $D([0, 1))$ be the open rigid analytic disk of radius 1 with coordinate $u$. The $\overline{K}$-points of $D([0, 1))$ correspond to elements $x \in \overline{K}$ such that $|x| < 1$. We can think of the points of $D([0,1))$ as the orbits of $\Gal(\overline{K}/K_{0})$ acting on the set $\{x  \in \overline{K} \text{ }|\text{ } |x|<1 \}.$ Let $\mathcal{O}:=\mathcal{O}_{D([0,1))}$ be the multivariable ring of rigid-analytic functions on the $D([0,1))$, which we can think of as a structure sheaves over $D([0, 1)^\delta)$.

Before introducing the category of $(\varphi, N_\nabla)$-modules over $\mathcal{O}$, we recall a fundamental ring-theoretic property of $\mathcal{O}$. Since $K_0$ is discretely valued, every invertible sheaf on the rigid-analytic open unit disc $\Delta$ is trivial, and every effective divisor on $\Delta$ is the divisor of an analytic function. Consequently, although $\mathcal{O}$ is not Noetherian, it is a Bézout domain, meaning that every finitely generated ideal of $\mathcal{O}$ is principal. This Bézout property ensures that coherent $\mathcal{O}$-modules behave much like modules over a principal ideal domain, as captured by the following key structural result (see \cite[Lemma~1.1.4]{Kis1}):

\begin{lemma}
    Let $I\subset [0,1)$ be an interval, $M$ a free $\mathcal{O}_I$-module of finite rank, and $N\subset M$ an $\mathcal{O}_I$-submodule.
Then the following three properties are equivalent:
    \begin{enumerate}
        \item $N\subset M$ is a closed submodule,
        \item $N$ is a finitely generated module, 
        \item $N$ is a free module.
    \end{enumerate}
\end{lemma}

\begin{dfn}
Let $\text{Mod}^{\varphi, N_{\nabla}}_{\mathcal{O}}$ be the category whose objects are finite free $\mathcal{O}$-modules $\mathcal{M}$ equipped with the following two operators:
\begin{enumerate}
    \item an $\varphi_{\mathcal{O}}$-semilinear, injective endomorphism $\varphi_{\mathcal{M}}$ of $\mathcal{M}$ such that the cokernel of $(1\otimes\varphi_{\mathcal{M}})$ is annihilated by a power $E^h$ for some integer $h\geq 0$.
    \item  a $K_0$-linear endomorphism $N_\nabla$ of $\mathcal{M}$ satisfying the Leibniz-rule and the relation $N_\nabla \circ \varphi = p \cdot \frac{E(u)}{E(0)} \cdot \varphi \circ N_{\nabla}$.
\end{enumerate}
The morphisms of this category are $\mathcal{O}$-module homomorphisms that are compatible with these additional structures. This $\text{Mod}^{\varphi, N_{\nabla}}_{\mathcal{O}}$ is equivalent the category of effective, filtered $(\varphi, N)$-modules over $K$, which was introduced by Fontaine (\cite{Fo2}).
\end{dfn}

The work of Zábrádi Gergely (\cite{Z4, Z5, Z6, Z7}) including joint work with Aprameyo Pal (\cite{PZ3}), Annie Carter and Kiran S. Kedlaya (\cite{CKZ1}), and Jishnu Ray and Feng Wei (\cite{RWZ2}), established the theory of multivariable $(\varphi, \Gamma)$-modules, which gives rise to the multivariable $p$-adic Hodge theory. These modules are used to study the $p$-adic representations of $G_{K, \Delta}:=\Gal(\overline{K}/K_{0})\times \cdots \times \Gal(\overline{K}/K_{0})$. 

In parallel with Zábrádi's work, the research of Scholze and Weinstein (\cite{SchWei1, Wei2}) helps us understand $G_{K, \Delta}$ as the fundamental group of the diamond $\text{Spd}\mathbb{Q}_p \times \cdots \times \text{Spd}\mathbb{Q}_p$. Building on these foundations, Olivier Brinon, Bruno Chiarellotto, and Nicola Mazzari (\cite{BCM1}) established the corresponding de Rham and Hodge-Tate period rings, multivariable Sen theory, and multivariable Hodge-Tate and de Rham $p$-adic representations of $G_{K, \Delta}$. Léo Poyeton and Pietro Vanni (\cite{PV1}) expanded upon this research by generalizing Colmez-Tate-Sen descent to the multivariable case, defining crystalline and semi-stable period rings, and introducing the concept of multivariable admissibility.

The most recent work in this area is by Aprameyo Pal and Rohit Pokhrel (\cite{Pal1}), who provide a cleaner way to study these objects using an approach closely related to the classical setting. Their article investigates the geometric side of the multivariable $p$-adic Hodge theory, although it does not build on Scholze's work in diamonds (\cite{SchWei1}). 

Our article aims to establish the rigid-analytic side of rational multivariable $p$-adic Hodge theory, generalizing Kisin's equivalence between effective filtered $(\varphi, N)$-modules and $(\varphi, N_\nabla)$-modules over multivariable rigid-analytic functions on the multivariable open rigid-analytic polydisk of radius $1$ over the tensor product of fields. In the next subsection, we will introduce the relevant notation and state the main theorem of the article.

\subsection*{Notation}
Throughout this paper, we use the following notation. Let $\Delta$ be an index set with cardinality $|\Delta| = \delta$. For each $i\in \Delta$, let $k_i$ be a perfect field characteristic $p$, let $W_i:=W(k_i)$ be the ring of Witt vectors, and let $K_{0, i}:= W_i\left[\frac{1}{p}\right]$ be the fraction field of $W_i$. Let $K_i/K_{0,i}$ be a finite, totally ramified field extension of $K_{0,i}$, and let  $\pi_{i}\in K_i$ be a uniformizer with $E_i(u_i) \in K_{0,i}[u_i]$ Eisenstein polynomial. Let $\overline{K_i}$ be an algebraic closure of $K_i$, and define a sequence of elements for $n \in \mathbb{Z}^{\geq 0}$ by setting $\pi_{0, i}:=\pi_i$, $\pi_{n+1, i}^p=\pi_{n, i}$. We denote $K_{n+1, i}:= K_i(\pi_{n,i})$. Setting $\mathfrak{S}_i:=W_i[[u_i]]$, we define $\widehat{\mathfrak{S}_{n,i}}$ as the completion of $K_{n+1, i}\otimes_{W_i} \mathfrak{S}_i$ at the maximal ideal $(u_i-\pi_{n,i})$. 

The preceding definitions and notation are the univariate cases of our multivariable setup. Let $K_{0,\Delta}$ be the tensor product of $K_{0,i}$ over $\mathbb{Q}_p$, let $K_\Delta$ be the tensor product of $K_i$ over $\mathbb{Q}_p$, and let $E_\Delta(u_1,\cdots, u_\delta)$ be the product of the polynomials $E_i(u_i)$. Furthermore, let $\mathfrak{S}_\Delta$ be the completed tensor product of $\mathfrak{S}_i$ over $\mathbb{Z}_p$, and let $\widehat{\mathfrak{S}_{\underline{n}, \Delta}}$ be the completed tensor product of $\widehat{\mathfrak{S}_{n_i, i}}$ over $\mathbb{Q}_p$ for  $(\underline{n}:=(n_1, \cdots, n_\delta))$. We can think of  $\mathfrak{S}_\Delta$ as $(W_1\otimes_{\mathbb{Z}_p}\cdots \otimes_{\mathbb{Z}_p} W_\delta)[[u_1, \cdots, u_\delta]]$.

We can explicitly describe these rings via their connected components, because each $W_i / \mathbb{Z}_p$ is a finite unramified extension, the base tensor product $W_\Delta := W_1 \otimes_{\mathbb{Z}_p} \dots \otimes_{\mathbb{Z}_p} W_\delta$ is a finite étale $\mathbb{Z}_p$-algebra, and therefore decomposes into a finite direct product of Witt rings $\prod_{\alpha \in \Sigma} W(k_\alpha)$ corresponding to its primitive idempotents $\{e_\alpha\}_{\alpha \in \Sigma}$. Consequently, $K_{0,\Delta} = W_\Delta[1/p] \cong \prod_{\alpha \in \Sigma} F_\alpha$, where each $F_\alpha := e_\alpha K_{0,\Delta}$ is a finite unramified extension of $\mathbb{Q}_p$. $\mathfrak{S}_\Delta \cong W_\Delta[[u_1, \dots, u_\delta]] \cong \prod_{\alpha \in \Sigma} W(k_\alpha)[[u_1, \dots, u_\delta]]$. Thus, $\mathfrak{S}_\Delta$ is a finite direct product of regular local rings of dimension $\delta + 1$. The partial Frobenii $\{\sigma_i\}_{i \in \Delta}$ act on $W_\Delta$ and $K_{0,\Delta}$, transitively permuting the idempotents $\{e_\alpha\}_{\alpha \in \Sigma}$. Any module over these rings similarly decomposes into components (e.g., $D = \bigoplus_{\alpha} e_\alpha D$). The next lemma will describe the structure of $K_{0, \Delta}$ and the connected components.

\begin{lemma}\label{lemma:components}
Assume that for each $i \in \Delta$, the residue field $k_i$ is a finite extension of $\mathbb{F}_p$ (so that $K_{0,i} = W(k_i)[1/p]$ is a finite unramified extension of $\mathbb{Q}_p$ with $K_{0,i}^{\sigma_i = 1} = \mathbb{Q}_p$), then:
\begin{enumerate}
    \item $K_{0,\Delta} = \bigotimes_{i \in \Delta, \mathbb{Q}_p} K_{0,i}$ is a finite étale $\mathbb{Q}_p$-algebra, which decomposes into a finite direct product of fields $K_{0,\Delta} \cong \prod_{\alpha \in \Sigma} F_\alpha$ corresponding to its primitive idempotents $\{e_\alpha\}_{\alpha \in \Sigma}$, where $F_\alpha = e_\alpha K_{0,\Delta}$.
    \item The partial Frobenii $\{\sigma_i\}_{i \in \Delta}$ permute the primitive idempotents $\{e_\alpha\}_{\alpha \in \Sigma}$ transitively, and $K_{0,\Delta}^{\sigma_i = 1 \;\forall i \in \Delta} = \mathbb{Q}_p$. Consequently, any finitely generated $K_{0,\Delta}$-module $D$ equipped with bijective $\sigma_i$-semilinear operators $\varphi_i$ ($i \in \Delta$) has constant dimension across all field factors $F_\alpha$, and is therefore a finite free $K_{0,\Delta}$-module.
\end{enumerate}
\end{lemma}

\begin{proof}
Since each $k_i / \mathbb{F}_p$ is a finite separable field extension, $K_{0,i} = W(k_i)[1/p]$ is a finite unramified (hence finite étale) field extension of $\mathbb{Q}_p$, because finite étale algebras over a field are stable under base change and tensor products, $K_{0,\Delta} = \bigotimes_{i \in \Delta, \mathbb{Q}_p} K_{0,i}$ is a finite étale $\mathbb{Q}_p$-algebra. By the structure theorem for finite étale algebras over a field, $K_{0,\Delta}$ is a finite-dimensional commutative semisimple $\mathbb{Q}_p$-algebra, and therefore decomposes uniquely into a finite direct product of finite field extensions of $\mathbb{Q}_p$:
\[
K_{0,\Delta} \cong \prod_{\alpha \in \Sigma} F_\alpha, \qquad F_\alpha := e_\alpha K_{0,\Delta},
\]
where $\{e_\alpha\}_{\alpha \in \Sigma}$ is the set of primitive orthogonal idempotents of $K_{0,\Delta}$ satisfying $e_\alpha^2 = e_\alpha$, $e_\alpha e_\beta = 0$ for $\alpha \neq \beta$, and $\sum_{\alpha \in \Sigma} e_\alpha = 1$.

For each $i \in \Delta$, the partial Frobenius $\sigma_i$ acts as the absolute Frobenius automorphism on $K_{0,i}$ and as the identity on $K_{0,j}$ for $j \neq i$, hence defines a $\mathbb{Q}_p$-algebra automorphism of $K_{0,\Delta}$. Any ring automorphism of $K_{0,\Delta}$ must permute its primitive idempotents $\{e_\alpha\}_{\alpha \in \Sigma}$. 

We first show by induction on $|\Delta| = \delta$ that the subring of simultaneously fixed elements $K_{0,\Delta}^{\sigma_i = 1 \;\forall i \in \Delta}$ is equal to $\mathbb{Q}_p$. For $\delta = 1$, we have $K_{0,1}^{\sigma_1 = 1} = W(k_1)^{\sigma_1 = 1}[1/p] = \mathbb{Z}_p[1/p] = \mathbb{Q}_p$. For $\delta > 1$, write $K_{0,\Delta} = K_{0,\delta} \otimes_{\mathbb{Q}_p} K_{0,\Delta \setminus \{\delta\}}$. Since $\sigma_\delta$ acts as $\sigma_\delta \otimes \mathrm{id}_{K_{0,\Delta \setminus \{\delta\}}}$ and $K_{0,\Delta \setminus \{\delta\}}$ is flat over $\mathbb{Q}_p$, taking $\sigma_\delta$-invariants commutes with tensoring over $\mathbb{Q}_p$:
\[
(K_{0,\Delta})^{\sigma_\delta = 1} = K_{0,\delta}^{\sigma_\delta = 1} \otimes_{\mathbb{Q}_p} K_{0,\Delta \setminus \{\delta\}} = \mathbb{Q}_p \otimes_{\mathbb{Q}_p} K_{0,\Delta \setminus \{\delta\}} = K_{0,\Delta \setminus \{\delta\}}.
\]
Taking invariants under the remaining automorphisms $\sigma_1, \dots, \sigma_{\delta-1}$ and applying the induction hypothesis yields $K_{0,\Delta}^{\sigma_i = 1 \;\forall i \in \Delta} = \mathbb{Q}_p$.

Now let $\mathcal{O} \subseteq \Sigma$ be any orbit of the action of the group $\langle \sigma_i \mid i \in \Delta \rangle$ on $\Sigma$, and consider the idempotent $e_{\mathcal{O}} := \sum_{\alpha \in \mathcal{O}} e_\alpha \in K_{0,\Delta}$. Since $\mathcal{O}$ is stable under $\sigma_i$ for every $i \in \Delta$, we have $\sigma_i(e_{\mathcal{O}}) = e_{\mathcal{O}}$ for all $i \in \Delta$, so $e_{\mathcal{O}} \in K_{0,\Delta}^{\sigma_i = 1 \;\forall i \in \Delta} = \mathbb{Q}_p$. Since the field $\mathbb{Q}_p$ contains only the trivial idempotents $0$ and $1$, and $e_{\mathcal{O}} \neq 0$, we must have $e_{\mathcal{O}} = 1$, which implies $\mathcal{O} = \Sigma$. Thus, the action of $\{\sigma_i\}_{i \in \Delta}$ on $\{e_\alpha\}_{\alpha \in \Sigma}$ is transitive.

Let $D$ be any finitely generated $K_{0,\Delta}$-module equipped with bijective $\sigma_i$-semilinear maps $\varphi_i : D \to D$. Using the idempotent decomposition $1 = \sum_{\alpha \in \Sigma} e_\alpha$, we have $D = \bigoplus_{\alpha \in \Sigma} e_\alpha D$, where each $e_\alpha D$ is a finite-dimensional vector space over the field $F_\alpha = e_\alpha K_{0,\Delta}$. For any $\alpha \in \Sigma$ and $i \in \Delta$, write $\sigma_i(e_\alpha) = e_{\sigma_i(\alpha)}$, then for any $x \in D$,
\[
\varphi_i(e_\alpha x) = \sigma_i(e_\alpha)\varphi_i(x) = e_{\sigma_i(\alpha)}\varphi_i(x) \in e_{\sigma_i(\alpha)} D,
\]
and since $\varphi_i$ is bijective and $\sigma_i|_{F_\alpha} : F_\alpha \xrightarrow{\sim} F_{\sigma_i(\alpha)}$ is a field isomorphism, $\varphi_i$ restricts to a $\sigma_i$-semilinear bijection $e_\alpha D \xrightarrow{\sim} e_{\sigma_i(\alpha)} D$. In particular, $\dim_{F_\alpha}(e_\alpha D) = \dim_{F_{\sigma_i(\alpha)}}(e_{\sigma_i(\alpha)} D)$ for all $i \in \Delta$. By the transitivity of the action of $\{\sigma_i\}_{i \in \Delta}$ on $\Sigma$, the dimension $d := \dim_{F_\alpha}(e_\alpha D)$ is independent of $\alpha \in \Sigma$. Choosing an $F_\alpha$-basis of $d$ elements in each $e_\alpha D$ and summing them over $\alpha \in \Sigma$ gives a $K_{0,\Delta}$-basis of $D$ of cardinality $d$, so $D \cong \prod_{\alpha \in \Sigma} F_\alpha^{\oplus d} \cong K_{0,\Delta}^{\oplus d}$ is a finite free $K_{0,\Delta}$-module of rank $d$.
\end{proof}

Let $D([0, 1)^\delta)$ be the multivariable open rigid-analytic polydisk of radius $1$ over $K_{0,\Delta}$ with coordinates $(u_1, \dots, u_\delta)$. Under the idempotent decomposition $K_{0,\Delta} \cong \prod_{\alpha \in \Sigma} F_\alpha$ of Lemma \ref{lemma:components}, $D([0, 1)^\delta)$ decomposes as a finite disjoint union $\coprod_{\alpha \in \Sigma} D([0, 1)^\delta)_{F_\alpha}$ of open unit polydisks over the fields $F_\alpha$. Set-theoretically, the geometric points of each component correspond to tuples $\mathbf{x} = (x_1, \dots, x_\delta) \in \prod_{i \in \Delta} \overline{K}_i$ such that $\vert{}x_i\vert{} < 1$ for all $i \in \Delta$, and the maximal ideals of the coordinate ring on $D([0, 1)^\delta)$ correspond to pairs consisting of a connected component $\alpha \in \Sigma$ and an orbit of $G_{K_0, \Delta} := \mathrm{Gal}(\overline{K}_1/K_{0,1}) \times \cdots \times \mathrm{Gal}(\overline{K}_\delta/K_{0,\delta})$ acting componentwise on $\{ \mathbf{x} \in \prod_{i \in \Delta} \overline{K}_i \mid \vert{}x_i\vert{} < 1 \text{ for all } i \in \Delta \}$.

Let $I_\Delta = (I_1, I_2, \dots, I_\delta)$ be a tuple of subintervals, where each $I_i \subset [0,1)$. The admissible open subspace $D(I_\Delta) \subset D([0,1)^\delta)$ is a poly-annulus over $K_{0,\Delta}$, whose geometric points correspond to tuples $\mathbf{x} = (x_1, \dots, x_\delta) \in \prod_{i \in \Delta} \overline{K}_i$ such that $|x_i| \in I_i$ for all $i \in \Delta$. Let $\mathcal{O}_{I_\Delta} := \Gamma(D(I_\Delta), \mathcal{O}_{D(I_\Delta)})$ be the Fr\'echet $K_{0,\Delta}$-algebra of rigid-analytic functions on the poly-annulus $D(I_\Delta)$. denoting by $\mathcal{O}_{I_i}$ the $K_{0,i}$-algebra of rigid-analytic functions on the univariate annulus $D(I_i)$ over $K_{0,i}$, we can view $\mathcal{O}_{I_\Delta}$ as the completed tensor product over $\mathbb{Q}_p$:
\[
\mathcal{O}_{I_\Delta} = \mathcal{O}_{I_1} \widehat{\otimes}_{\mathbb{Q}_p} \cdots \widehat{\otimes}_{\mathbb{Q}_p} \mathcal{O}_{I_\delta} \cong (\mathcal{O}_{I_1} \otimes_{K_{0,1}} K_{0,\Delta}) \widehat{\otimes}_{K_{0,\Delta}} \cdots \widehat{\otimes}_{K_{0,\Delta}} (\mathcal{O}_{I_\delta} \otimes_{K_{0,\delta}} K_{0,\Delta}).
\]
Setting $\mathcal{O}_\Delta := \mathcal{O}_{[0,1)^\delta}$ and $\mathcal{O}_i := \mathcal{O}_{[0,1)_i}$, we can view $\mathcal{O}_{I_\Delta}$ and $\mathcal{O}_\Delta$ as the rings of global sections of the structure sheaves on $D(I_\Delta)$ and $D([0, 1)^\delta)$, respectively.

\begin{lemma}\label{lemma:components2}
    For every closed affinoid polydisk $U_{\rho_k} \subset D([0,1)^\delta)$, the ring $\mathcal{O}_\Delta(U_{\rho_k}) \cong \prod_{\alpha \in \Sigma} e_\alpha \mathcal{O}_\Delta(U_{\rho_k})$ is a finite direct product of free Tate algebras over the fields $F_\alpha$ (hence a finite direct product of regular Noetherian normal domains), and $\mathcal{O}_\Delta \cong \prod_{\alpha \in \Sigma} e_\alpha \mathcal{O}_\Delta$.
\end{lemma}

\begin{proof}
Since $1 = \sum_{\alpha \in \Sigma} e_\alpha$ is a finite orthogonal idempotent decomposition in the base ring $K_{0,\Delta} \subset \mathcal{O}_\Delta(U_{\rho_k})$, any $\mathcal{O}_\Delta(U_{\rho_k})$-algebra or module splits as the finite direct product of its $e_\alpha$-components:
\[
\mathcal{O}_\Delta(U_{\rho_k}) = \prod_{\alpha \in \Sigma} e_\alpha \mathcal{O}_\Delta(U_{\rho_k}).
\]
For each $\alpha \in \Sigma$, $e_\alpha \mathcal{O}_\Delta(U_{\rho_k})$ is the ring of rigid-analytic functions on the closed affinoid polydisk $U_{\rho_k}$ over the $p$-adic field $F_\alpha = e_\alpha K_{0,\Delta}$, which is isomorphic to a $\delta$-dimensional free Tate algebra over the field $F_\alpha$ \cite[Section 6.1]{BGR1}. Since every free Tate algebra over a non-Archimedean field is a regular Noetherian UFD (and hence a Noetherian normal domain) \cite[Section 5.2.6, Theorem 1]{BGR1}, $\mathcal{O}_\Delta(U_{\rho_k})$ is a finite direct product of regular Noetherian normal domains. Taking the projective limit over $k \to \infty$ commutes with the finite direct product over $\Sigma$, yielding $\mathcal{O}_\Delta = \varprojlim_k \mathcal{O}_\Delta(U_{\rho_k}) \cong \prod_{\alpha \in \Sigma} e_\alpha \mathcal{O}_\Delta$.
\end{proof}

The main purpose of this article is to introduce a generalization of $(\varphi, N)$-modules and $(\varphi, N_\nabla)$-modules. The following definitions of these two categories are not entirely precise, but they will allow us to state the main theorem of the article.

\vspace{0.5cm}

\begin{dfn}
 Let $\text{MF}_{K_{0, \Delta}}^{(\varphi_i, N_i)_{i\in \Delta}, \text{ Fil}\geq 0}$ be the category of effective, filtered $(\varphi_i, N_i)_{i\in \Delta}$-modules over $K_{0,\Delta}$. The objects of this category are finite free $K_{0,\Delta}$-modules $D$ equipped with two sets of operators, $\{\varphi_i\}_{i\in\Delta}$ and $\{N_i\}_{i\in \Delta}$. These operators satisfy the same properties as in the classical case and commute with one another for all indices, with the sole exception that $N_i\circ \varphi_i = p\cdot(\varphi_i\circ N_i)$ for each $i\in \Delta$. Furthermore, the module $D_{K_\Delta}:=D\otimes_{K_{0, \Delta}}K_\Delta$ is equipped with a decreasing, separated, and exhaustive multifiltration. The morphisms between the objects are strictly filtered $K_0$-linear homomorphisms which commute with $\varphi_i$ and $N_i$ for every $i\in \Delta$.
\end{dfn}

\vspace{0.5cm}

\begin{dfn}
 Let $\text{Mod}_{\mathcal{O}_\Delta}^{(\varphi_i, N_\nabla^{i})_{i\in \Delta}}$ be the category of $(\varphi_i, N_\nabla^i)_{i\in \Delta}$-modules over $\mathcal{O}_\Delta$. The objects of this category are finite free $\mathcal{O}_\Delta$-modules $\mathcal{M}$ equipped with two sets of operators, $\{\varphi_i\}_{i\in\Delta}$ and $\{N_\nabla^i\}_{i\in \Delta}$. These operators satisfy the same properties as in the classical case and commute with one another for all indices, with the sole exception that $N_\nabla^i\circ \varphi_i = p\cdot\frac{E_i(u_i)}{E_i(0)}\cdot(\varphi_i\circ N_\nabla^i)$ for each $i\in \Delta$. Furthermore, we define a total Frobenius map on $M$ by $\varphi_\Delta:=\varphi_1\circ\varphi_2\circ\ldots\circ\varphi_\delta$. We require that there exists an $\underline{h}:=(h_1, \cdots, h_\delta)\in \mathbb{Z}^\delta$ such that the cokernel of the linearizad map $1\otimes \varphi_\Delta$ is annihilated by the ideal $(E_1^{h_1}(u_1), \dots,  E_\delta^{h_\delta}(u_\delta))$, meaning 
 $(E_1^{h_1}(u_1), \dots,  E_\delta^{h_\delta}(u_\delta))\cdot \text{coker}(1\otimes\varphi_\Delta)=0.$ The morphisms of this category are $\mathcal{O}_\Delta$-linear homomorphisms that commute with $\varphi_i$ and $N_\nabla^i$ for every $i\in \Delta$.
\end{dfn}

\vspace{0.5cm}

We can define the functor from category $\text{Mod}_{\mathcal{O}_\Delta}^{(\varphi_i, N_\nabla^{i})_{i\in \Delta}}$ to category $\text{MF}_{K_{0, \Delta}}^{(\varphi_i, N_i)_{i\in \Delta}, \text{ Fil}\geq 0}$ by taking the reduction modulo the ideal $(u_1, \cdots, u_\delta)$, and we denote this functor by $\mathcal{D}$. The functor in the opposite direction, denoted by $\mathcal{M}$, associates to a multivariable filtered $(\varphi_i, N_i)$-module $D$ a multivariable $(\varphi_i, N_\nabla^i)$-module over $\mathcal{O}_{\Delta}$, constructed as the intersection of the base-changed univariate $(\varphi_i, N_\nabla^i)$-modules $\mathcal{M}_i(D)$ along each variable direction.

\vspace{0.5cm}

\begin{thm}
    The categories $\text{MF}_{K_{0,\Delta}}^{(\varphi_i, N_i)_{i\in \Delta}, \text{ Fil}\geq 0}$ and $\text{Mod}_{\mathcal{O}_\Delta}^{(\varphi_i, N_\nabla^{i})_{i\in \Delta}}$ are exact, tensor categories. The functors $\mathcal{D}$ and $\mathcal{M}$ are quasi-inverse to each other, establishing an equivalence of categories.
\end{thm}

\vspace{1.0cm}
\noindent\textbf{Acknowledgements.} I would like to thank my supervisor Gergely Zábrádi for suggesting this problem, for his guidance throughout its solution, and for his careful reading of all previous versions of the manuscript.
\vspace{1.0cm}

\noindent\textbf{AI declaration.} The authors acknowledge the use of Google’s Gemini models in developing some of the arguments, accessed through a student subscription. The key elements in the proofs are the authors’ own. The authors take sole responsibility for any mistakes in this work.

\section*{Preliminaries}
\subsection*{Reflexive modules}
We first recall some properties of reflexive modules that will be used in the following sections. In linear algebra, every finite-dimensional vector space $V$ is naturally isomorphic to its double dual. However, when considering modules over a general ring, the double dual loses information about the torsion, and the natural map from a module to its double dual is not necessarily surjective. Modules for which this natural map is an isomorphism are called reflexive modules. We will demonstrate why this property is useful, drawing on Commutative Algebra book of Bourbaki (\cite{B1}) and the Stacks project (\cite{stacks-project}).

\begin{dfn}
Let $R$ be a domain. An $R$-module $M$ is reflexive if the natural map $j:M \rightarrow\text{Hom}_R(\text{Hom}_R(M,R),R)$, which sends $m\in M$ to the evaluation map $\phi \mapsto \phi(m)$ for all $\phi\in \text{Hom}_R(M,R)$, is an isomorphism. Furthermore, we denote $\text{Hom}_R(M,R)$ by $M^*$, and we refer to  $\text{Hom}_R(\text{Hom}_R(M,R),R)$ as the reflexive hull of $M$. (\cite[Definition 15.24.1]{stacks-project})
\end{dfn}

The definition and the brief introduction give us some hint of the statement of the following lemma.

\begin{lemma}
Let $R$ be a domain and let $M$ be an $R$-module, then the following are true:
\begin{enumerate}
    \item If $M$ reflexive, then $M$ is torsion-free.
    \item If $M$ is finite, then the kernel and cokernel of $j$ map are torsion modules.
    \item If $M$ is finite, then $j$ is injective if and only if $M$ is torsion free.
\end{enumerate}
Furthermore, if $R$ is a DVR and $M$ finite, then $j$ will be surjective. (\cite[Lemma 15.24.2]{stacks-project})
\end{lemma}

The following lemma shows that local reflexiveness implies global reflexiveness.

\begin{lemma}
Let $R$ be a Noetherian domain and let $M$ be a finite $R$-module, then the following are equivalent:
\begin{enumerate}
    \item $M$ is reflexive.
    \item $M_{\mathfrak{p}}$ is reflexive $R_{\mathfrak{p}}$-module for all prime ideals of $R$.
    \item $M_{\mathfrak{m}}$ is a reflexive $R_{\mathfrak{m}}$-module for all maximal ideals of $R$. 
\end{enumerate}
(\cite[Lemma 15.24.4]{stacks-project})
\end{lemma}

\begin{dfn}
A finitely generated module $M$ over a Noetherian ring $R$ satisfies the $S_2$ condition if, for every prime ideal $\mathfrak{p}$ in $R$, the algebraic depth of the localization $M_{\mathfrak{p}}$ satisfies:$$\text{depth}(M_{\mathfrak{p}}) \ge \min(2, \text{dim}(M_{\mathfrak{p}})).$$ 
\end{dfn}

This condition can be generalized to any $k \in \mathbb{N}$, and because it was introduced by Jean-Pierre Serre, it is commonly known as Serre's condition $S_k$. These conditions give a necessary and sufficient conditions for a commutative Noetherian ring to be normal, provided we also require that for any prime ideal of height at most $k$, the corresponding localization is a regular local ring.

\begin{lemma}
Let $R$ be a Noetherian domain, and let $M$ be a finite $R$-module, then the following are equivalent:
\begin{enumerate}
    \item $M$ is reflexive.
    \item For every $\mathfrak{p}$ of $R$, either $M_{\mathfrak{p}}$ is a reflexive $R_{\mathfrak{p}}$-module, or $\text{depth}(M_{\mathfrak{p}})\geq 2$.
\end{enumerate}
(\cite[Lemma 15.24.15]{stacks-project})
\end{lemma}

\begin{lemma}
Let $R$ be a Noetherian domain, and let $M$ be a finite reflexive $R$-module. Let $\mathfrak{p} \subset R$ be a prime ideal. 
\begin{enumerate}
    \item If depth$(R_{\mathfrak{p}})\geq 2$, then depth$(M_{\mathfrak{p}})\geq 2$.
    \item If $R$ satisfies $(S_2)$, then $M$ satisfies $(S_2)$.
\end{enumerate}
(\cite[Lemma 15.24.16]{stacks-project})
\end{lemma}

This previous lemmas and the following lemma show that there is some connection between $(S_2)$ condition and reflexivity. However, the following lemma demonstrates that this implication only holds in one direction, \cite[Example 15.24.17]{stacks-project} provides a counterexample for the converse.

\begin{lemma}\label{lemma:reflexiv}
    Let $R$ be a Noetherian normal domain with fraction field $K$, and $M$ be a finitely generated $R$-module, then following are equivalent:
    \begin{enumerate}
        \item $M$ is a reflexive module,
        \item $M$ is torsion free and has property $(S_2)$,
        \item $M$ is torsion-free and 
        \[
        M = \bigcap_{\operatorname{ht}(\mathfrak{p})=1} M_{\mathfrak{p}},
        \]
        where the intersection is taken inside $M_K=M\otimes_{R} K$ over all prime ideals $\mathfrak{p}$ of height $1$.
    \end{enumerate}
(\cite[Lemma 15.24.18]{stacks-project})
\end{lemma}

More generally, if $R = \prod_{\alpha \in \Sigma} R_\alpha$ is a finite direct product of Noetherian normal domains with total ring of fractions $Q(R) := \prod_{\alpha \in \Sigma} \mathrm{Frac}(R_\alpha)$ (such as $R = \mathcal{O}_\Delta(U_{\rho_k})$ by Lemma \ref{lemma:components}), then every $R$-module $M$ decomposes canonically as $M = \bigoplus_{\alpha \in \Sigma} e_\alpha M$ where $e_\alpha M$ is an $R_\alpha$-module, and every height-1 prime ideal $\mathfrak{p}$ of $R$ corresponds to a height-1 prime ideal $\mathfrak{p}_\alpha$ of a unique factor $R_\alpha$ with $M_{\mathfrak{p}} \cong (e_\alpha M)_{\mathfrak{p}_\alpha}$. Applying Lemma \ref{lemma:reflexiv} componentwise to each factor $R_\alpha$ shows that a finitely generated torsion-free $R$-module $M$ is reflexive if and only if each $e_\alpha M$ is reflexive over $R_\alpha$, which holds if and only if
\[
M = \bigoplus_{\alpha \in \Sigma} \bigcap_{\substack{\mathfrak{p}_\alpha \subset R_\alpha \\ \mathrm{ht}(\mathfrak{p}_\alpha) = 1}} (e_\alpha M)_{\mathfrak{p}_\alpha} = \bigcap_{\mathrm{ht}(\mathfrak{p}) = 1} M_{\mathfrak{p}} \subset M \otimes_R Q(R).
\]
In particular, if $M_1, M_2 \subset M \otimes_R Q(R)$ are two finitely generated reflexive (for instance, free) $R$-submodules of the same generic space, then $M_1 = M_2$ if and only if $(M_1)_{\mathfrak{p}} = (M_2)_{\mathfrak{p}}$ for every height-1 prime ideal $\mathfrak{p}$ of $R$.

In our multivariable setting, this last argument provides the essential local tool for studying intersections of partial $(\varphi_i, N_\nabla^i)$-modules: it allows us to verify local freeness and equality of lattices over Noetherian normal affinoid subdomains simply by checking their localizations at height-1 prime ideals. However, because our multivariable $(\varphi_i, N_\nabla^i)$-module over $\mathcal{O}_{\Delta}$ live over the non-Noetherian Fréchet algebra $\mathcal{O}_{\Delta}$ on the open polydisk $D([0,1)^{\delta})$, purely algebraic methods over affinoid rings are not sufficient on their own. To pass from these local affinoid properties to a globally free $\mathcal{O}_{\Delta}$-module, we need the rigid-analytic machinery of quasi-Stein spaces and Kiehl's Theorems A and B.

\subsection*{Kiehl's Theorem A and Theorem B}

For each $i\in \Delta$, let $I_i = (a_i, b_i) \subset [0,1)$ be an open interval, and $\{J_{i,n} = [r_{i,n}, s_{i,n}]\}_{n=1}^\infty$ be a sequence of nested closed intervals in $I_i$ such that $J_{i,n} \subset J_{i,n+1}$ and $\bigcup_{n=1}^\infty J_{i,n} = I_i$. For each $n\geq 1$, we define a closed affinoid poly-annulus $D_n := \{ \mathbf{x} \in D([0,1)^\delta) \mid \vert{}x_i\vert{} \in J_{i,n} \text{ for all } i \in \Delta \}$. This  $D(I_\Delta)$, where $I_\Delta = (I_1, I_2, \dots, I_\delta)$ is the admissible union of this nested ascending chain of affinoid spaces. 

The following norms make the $\mathcal{O}_{I_\Delta}$ into a Fréchet space. For any $f\in \mathcal{O}_{I_\Delta}$, we define the supremum norm on $D_n$ by $\Vert{}f\Vert{}_n := \sup_{\mathbf{x} \in D_n} \vert{}f(\mathbf{x})\vert{} < \infty$. The topology of $\mathcal{O}_{I_\Delta}$ is defined by uniform convergence on the affinoid subdomains $D_n$ as $n\rightarrow \infty$, because $\mathcal{O}_{I_\Delta}$ is complete with respect to this countable family of separating norms $\{ \Vert{} \cdot \Vert{}_n \}_{n=1}^\infty$, it forms a Fréchet algebra.

To work effectively with modules over our multivariable ring of functions $\mathcal{O}_{I_\Delta}$, we must rely on the foundational framework of rigid-analytic geometry. While the open disk $D([0,1)^\delta)$ and the poly-annulus $D(I_\Delta)$ are constructed via specific completed tensor products, they possess excellent geometric properties. 

\begin{dfn}
A quasi-Stein space is an adic space $X$ such that $\mathcal{O}(X)$ contains a topologically nilpotent unit, and $X$ can be written as the union of an ascending sequence $U_0 \subseteq U_1 \subseteq \cdots$ of adic affinoid subspaces where, for each $i\geq 0$, the transition map $\mathcal{O}(U_{i+1}) \rightarrow \mathcal{O}(U_i)$ has dense image with respect to the supremum norm on $U_k$. Note that by quasi-compactness, every adic affinoid subspace of $X$ is contained in some $U_i$. Moreover, $\mathcal{O}(X)= \varprojlim_{i} \mathcal{O}(U_i)$ is an open mapping ring.
\end{dfn}

Standard examples of quasi-Stein spaces include affinoid spaces, open rigid-analytic unit disks, one-dimensional annuli, and finite products of quasi-Stein spaces. In particular, both the open polydisk $D([0,1)^\delta)$ and the poly-annulus $D(I_\Delta)$ (exhausted by the nested closed affinoid subdomains $U_{\rho_k}$ and $D_n$, respectively, where polynomial functions are dense in each affinoid algebra) are quasi-Stein spaces, and their rings of global sections $\mathcal{O}_\Delta = \varprojlim_k \mathcal{O}_\Delta(U_{\rho_k})$ and $\mathcal{O}_{I_\Delta} = \varprojlim_n \mathcal{O}_{I_\Delta}(D_n)$ are Fréchet algebras.

\begin{dfn}
Let $(X, \mathcal{O}_X)$ be a ringed space. A sheaf of $\mathcal{O}_X$-modules $\mathcal{F}$ is said to be quasi-coherent if every point in $X$ admits an open neighborhood $U$ over which there exists an exact sequence
\[
\mathcal{O}_X^{\oplus I}|_U \longrightarrow \mathcal{O}_X^{\oplus J}|_U \longrightarrow \mathcal{F}|_U \longrightarrow 0
\]
for some index sets $I$ and $J$.
\end{dfn}

\begin{dfn}
Let $(X, \mathcal{O}_X)$ be a ringed space. A sheaf of $\mathcal{O}_X$-modules $\mathcal{F}$ is said to be coherent if it satisfying the following two properties:
\begin{enumerate}
    \item $\mathcal{F}$ is finite type over $\mathcal{O}_X$, that is, every point in $X$ has an open neighborhood $U$ in $X$ such that there is a surjective morphism $\mathcal{O}_X^n|_U \longrightarrow \mathcal{F}|_U$ for some $n\in\mathbb{N}$.
    \item For any open set $U\subset X$, any $n\in \mathbb{N}$, and any morphism $\alpha: \mathcal{O}_{X}^n|_U \longrightarrow \mathcal{F}|_U$ of $\mathcal{O}_X$-modules, the kernel of $\alpha$ is finite type.
\end{enumerate}
\end{dfn}

Based on the previous definition our $\mathcal{O}_\Delta$-modules are coherent sheaves over $D([0,1)^\delta)$, so we can use the following theorem by Kiehl (Theorem A, B) which is the rigid-analytic analogues of Cartan's Theorems A and B to pass from local affinoid data to global sections. This theorem provide the exact machinery required for understanding global structures on these spaces by analyzing local data and gluing it together, because it guarantee that the resulting global coherent sheaves behave predictably, possessing globally generated stalks and vanishing higher cohomology.

\begin{thm}[Kiehl's Theorem A and Theorem B]\label{thm:kiehl}
Let $X$ be a quasi-Stein space, and let $\mathcal{F}$ be a coherent sheaf of $\mathcal{O}_X$-module on $X$. Let $\{U_i\}$ be an admissible affinoid covering of $X$, then the following are true:
\begin{enumerate}
    \item For any $U_i$ in the admissible affinoid covering, the restriction map $\mathcal{F}(X) \rightarrow \mathcal{F}(U)$ has dense image.
    \item (Theorem A) For each $x\in X$, the stalk $\mathcal{F}_x$ is generated as an $\mathcal{O}_{X,x}$-module by global sections of $\mathcal{F}$.
    \item (Theorem B) For all $j>0$, the groups $H^j(X, \mathcal{F})=0$.
\end{enumerate}
(\cite{Kie1} and \cite{Kie2})
\end{thm}

When working with coherent sheaves over a quasi-Stein space $X = \bigcup_{k \ge 1} U_k$, Kiehl's Theorem B guarantees the vanishing of higher cohomology, but establishing that a coherent sheaf is globally free requires aligning local bases on the overlaps of the affinoid covering. This alignment is achieved via the following standard matrix factorization tool, originally due to Cartan and adapted to rigid analytic geometry.

\begin{lemma}[Cartan's Matrix Factorization]\label{lemma:cartan_factorization}
Let $V$ be a rigid analytic space, and let $U \subset V$ be a strictly inner affinoid subdomain with coordinate ring $\mathcal{O}(U)$ and supremum norm $\|\cdot\|_U$. For any integer $d \ge 1$, there exists a constant $\epsilon > 0$ such that any matrix $E \in \mathrm{GL}_d(\mathcal{O}(U))$ satisfying 
\[
\|E - I_d\|_U \le \epsilon
\]
admits a factorization
\[
E = H_U^{-1} H_V,
\]
where $H_U \in \mathrm{GL}_d(\mathcal{O}(U))$ and $H_V \in \mathrm{GL}_d(\mathcal{O}(V))$.
\end{lemma}

By combining this factorization with the density of $\mathcal{O}(V)$ in $\mathcal{O}(U)$, any transition matrix on the overlap of nested affinoids can be iteratively factored, allowing local bases of a projective module to be glued into a global basis over the entire quasi-Stein space (see \cite[Section 2.7]{FvdP1} or \cite[Section 6.1]{BGR1}).

If $X=\text{Sp}(A)$ is an affinoid space, then by Theorem A, every coherent sheaf of $\mathcal{O}_X$-modules is associated to a finitely generated $A$-module. This yields an equivalence between the category of holomorphic vector bundles on $X$ and the category of finitely generated projective $A$-modules. 

Let $X$ and $Y$ be holomorphic spaces, and let $\mathcal{F}$ be a vector bundle on $X\times Y$, then $\mathcal{F}$ is said to be locally trivial over $X$ if there exists a holomorphic covering $\mathfrak{U}=\{X_i\}_i$ of $X$ such that the restrictions $\mathcal{F}\vert{}_{X_i\times Y}$ are trivial..
\begin{thm}
Every holomorphic vector bundle on $X \times E^n$ is locally trivial over $X$.
(\cite[Theorem 1.]{Ltk1})
\end{thm}

Specializing this theorem to the case where $X = \mathrm{Sp}(K_0)$ is a single point, any open covering of $X$ is trivial, which implies that every holomorphic vector bundle on the $n$-dimensional closed unit polydisk $E^n = \mathrm{Sp}(T_n)$ is globally trivial. Combining this with the categorical equivalence between holomorphic vector bundles on $\mathrm{Sp}(T_n)$ and finitely generated projective $T_n$-modules established by Kiehl's Theorem A, we immediately obtain the rigid-analytic analogue of the Quillen--Suslin theorem:

\begin{cor}[Quillen-Suslin over Tate algebras]\label{cor:quillen-suslin}
Every finitely generated, projective modules over $n$ dimensional Tate algebra are free. (\cite[Theorem 1. Corollary]{Ltk1})
\end{cor}

In our multivariable setting over $K_{0,\Delta}$, a closed affinoid polydisk $U_{\rho_k} \subset D([0,1)^\delta)$ is not connected when $K_{0,\Delta}$ is not a field, as its coordinate ring decomposes into a finite direct product of free Tate algebras over the field factors of $K_{0,\Delta}$ by Lemma \ref{lemma:components2}. Applying Corollary \ref{cor:quillen-suslin} componentwise yields the following generalization:

\begin{cor}\label{cor:lutkebohmert_multivariable}
Let $U_{\rho_k} \subset D([0,1)^\delta)$ be a closed affinoid polydisk over $K_{0,\Delta}$, and let $P$ be a finitely generated projective $\mathcal{O}_\Delta(U_{\rho_k})$-module. If $P$ has constant rank $d$ across all connected components of $U_{\rho_k}$ (meaning that $\mathrm{rk}_{e_\alpha \mathcal{O}_\Delta(U_{\rho_k})}(e_\alpha P) = d$ for every primitive idempotent $e_\alpha \in K_{0,\Delta}$), then $P$ is a globally free $\mathcal{O}_\Delta(U_{\rho_k})$-module of rank $d$.
\end{cor}

\begin{proof}
Set $A := \mathcal{O}_\Delta(U_{\rho_k})$. By Lemma \ref{lemma:components2}, the primitive orthogonal idempotents $\{e_\alpha\}_{\alpha \in \Sigma}$ of $K_{0,\Delta} \cong \prod_{\alpha \in \Sigma} F_\alpha$ decompose $A$ into a finite direct product $A = \prod_{\alpha \in \Sigma} A_\alpha$, where each $A_\alpha := e_\alpha A$ is isomorphic to a $\delta$-dimensional free Tate algebra over the field $F_\alpha$. 

Since $1 = \sum_{\alpha \in \Sigma} e_\alpha$, the projective $A$-module $P$ decomposes canonically as $P = \bigoplus_{\alpha \in \Sigma} P_\alpha$, where each $P_\alpha := e_\alpha P$ is a finitely generated projective $A_\alpha$-module of rank $d$. Applying Corollary \ref{cor:quillen-suslin} to the free Tate algebra $A_\alpha$ over the field $F_\alpha$ for each $\alpha \in \Sigma$, we deduce that $P_\alpha$ is a free $A_\alpha$-module of rank $d$. 

For each $\alpha \in \Sigma$, choose an $A_\alpha$-basis $\{b_{\alpha, 1}, \dots, b_{\alpha, d}\}$ of $P_\alpha \subset P$, and define the global elements
\[
b_j := \sum_{\alpha \in \Sigma} b_{\alpha, j} \in P \qquad \text{for } j = 1, \dots, d.
\]
Since the idempotents $\{e_\alpha\}_{\alpha \in \Sigma}$ are mutually orthogonal and $b_{\alpha, j} \in e_\alpha P$, for any coefficients $c_j = \sum_{\alpha \in \Sigma} c_{\alpha, j} \in \prod_{\alpha \in \Sigma} A_\alpha = A$ (with $c_{\alpha, j} \in A_\alpha$) we have
\[
\sum_{j=1}^d c_j b_j = \sum_{\alpha \in \Sigma} \sum_{j=1}^d c_{\alpha, j} b_{\alpha, j}.
\]
Because $\{b_{\alpha, 1}, \dots, b_{\alpha, d}\}$ is an $A_\alpha$-basis of $P_\alpha$ for each $\alpha \in \Sigma$, it follows immediately that $\{b_1, \dots, b_d\}$ is an $A$-basis of $P = \bigoplus_{\alpha \in \Sigma} P_\alpha$, so $P \cong A^{\oplus d}$ is a free $\mathcal{O}_\Delta(U_{\rho_k})$-module of rank $d$.
\end{proof}

\subsection*{Multifiltrations}

To study the category of multivariable filtered $(\varphi_i, N_i)_{i \in \Delta}$-modules, we need a solid algebraic framework for filtrations indexed by partially ordered abelian groups. Following the formalism of Burgos Gil and Mallick \cite{GM1}, we briefly recall the foundational properties of multifiltrations specialized to the free abelian group $\mathbb{Z}^\delta$ (where $\delta = |\Delta|$) equipped with the standard strictly positive unimodular cone
\[
C := (\mathbb{Z}_{\ge 0})^\delta = \{ \underline{k} = (k_1, \dots, k_\delta) \in \mathbb{Z}^\delta \mid k_i \ge 0 \text{ for all } i \in \Delta \}.
\]
The cone $C$ induces a compatible partial order on $\mathbb{Z}^\delta$ defined by $\underline{l} \le \underline{k}$ if and only if $\underline{k} - \underline{l} \in C$, which means $l_i \le k_i$ for every $i \in \Delta$. Note that $(\mathbb{Z}^\delta, \le)$ forms a lattice in the poset sense, with component-wise meet $\underline{k} \wedge \underline{l} = (\min(k_i, l_i))_{i \in \Delta}$ and join $\underline{k} \vee \underline{l} = (\max(k_i, l_i))_{i \in \Delta}$.

For any subset $S \subset \mathbb{Z}^\delta$, we denote its upward closure with respect to $\le$ by
\[
\widetilde{S} := \{ \underline{k} \in \mathbb{Z}^\delta \mid \underline{k} \ge \underline{l} \text{ for some } \underline{l} \in S \}.
\]

\begin{dfn}
Let $E$ be a finitely generated, free module over $K_\Delta$ (such as $D_{K_\Delta} := D \otimes_{K_{0,\Delta}} K_\Delta$). A decreasing multifiltration $\mathrm{Fil}^{\bullet}$ on $E$ indexed by $(\mathbb{Z}^\delta, C)$ is a family of $K_\Delta$-submodules $\{\mathrm{Fil}^{\underline{k}}(E)\}_{\underline{k} \in \mathbb{Z}^\delta}$ such that whenever $\underline{l} \le \underline{k}$ in $\mathbb{Z}^\delta$, we have
\[
\mathrm{Fil}^{\underline{k}}(E) \subseteq \mathrm{Fil}^{\underline{l}}(E).
\]
For any subset $S \subset \mathbb{Z}^\delta$, we set $\mathrm{Fil}^S(E) := \sum_{\underline{k} \in S} \mathrm{Fil}^{\underline{k}}(E) = \mathrm{Fil}^{\widetilde{S}}(E)$.
\end{dfn}

While the classical notion of exhaustivity carries over directly from one-variable filtrations, separatedness and the structure of the associated graded object require extra care in the multivariable setting \cite{GM1}.

\begin{dfn}
Let $\mathrm{Fil}^{\bullet}$ be a decreasing multifiltration on $E$ indexed by $(\mathbb{Z}^\delta, C)$.
\begin{enumerate}
    \item $\mathrm{Fil}^{\bullet}$ is called exhaustive if $E = \sum_{\underline{k} \in \mathbb{Z}^\delta} \mathrm{Fil}^{\underline{k}}(E)$ (or equivalently, since $C$ is generating, $E = \bigcup_{\underline{k} \in \mathbb{Z}^\delta} \mathrm{Fil}^{\underline{k}}(E)$).
    \item $\mathrm{Fil}^{\bullet}$ is called chain-separated if for every strictly increasing chain $\underline{k}^{(1)} < \underline{k}^{(2)} < \dots$ in $\mathbb{Z}^\delta$, there exists an index $m_0$ such that $\mathrm{Fil}^{\underline{k}^{(m)}}(E) = 0$ for all $m \ge m_0$.
    \item $\mathrm{Fil}^{\bullet}$ is called separated if for every collection $\{S_\alpha\}$ of subsets of $\mathbb{Z}^\delta$ satisfying $\bigcap_\alpha \widetilde{S}_\alpha = \emptyset$, we have
    \[
    \bigcap_\alpha \mathrm{Fil}^{S_\alpha}(E) = \{0\}.
    \]
    \item $\mathrm{Fil}^{\bullet}$ is called regular if for every pair of subsets $S_1, S_2 \subset \mathbb{Z}^\delta$, the following distributivity condition holds:
    \[
    \left( \sum_{\underline{k} \in S_1} \mathrm{Fil}^{\underline{k}}(E) \right) \cap \left( \sum_{\underline{l} \in S_2} \mathrm{Fil}^{\underline{l}}(E) \right) = \sum_{\underline{m} \in \widetilde{S}_1 \cap \widetilde{S}_2} \mathrm{Fil}^{\underline{m}}(E).
    \]
    \item We say that $\mathrm{Fil}^{\bullet}$ is effective (denoted $\mathrm{Fil}^{\ge \underline{0}}$) if $\mathrm{Fil}^{\underline{0}}(E) = E$, where $\underline{0} = (0, \dots, 0) \in \mathbb{Z}^\delta$.
\end{enumerate}
\end{dfn}

Since $C = (\mathbb{Z}_{\ge 0})^\delta$ is a strongly strict polyhedral generating cone, any separated multifiltration on a finite-rank module $E$ is automatically chain-separated and admits only finitely many jump points \cite[Lemma 3.19 and Proposition 3.35]{GM1}. Conversely, any chain-separated multifiltration with finitely many jump points is separated \cite[Proposition 3.37]{GM1}.

In our multivariable setting, the ring $K_\Delta = \bigotimes_{i \in \Delta, \mathbb{Q}_p} K_i$ is a finite étale $K_{0,\Delta}$-algebra, and hence decomposes into a finite direct product of fields $K_\Delta \cong \prod_{\beta \in \Sigma'} L_\beta$ corresponding to its primitive idempotents $\{e'_\beta\}_{\beta \in \Sigma'}$. Applying the results of Burgos Gil and Mallick \cite{GM1} componentwise over each field factor $L_\beta$, the regularity condition prevents the pathological behavior of incomparable filtration jumps:

\begin{thm}\label{thm:multigrading}
Let $E$ be a finitely generated, free $K_\Delta$-module equipped with a decreasing multifiltration $\mathrm{Fil}^{\bullet}$ indexed by $(\mathbb{Z}^\delta, C)$, then $\mathrm{Fil}^{\bullet}$ is exhaustive, separated, and regular if and only if there exists a (non-canonical) multigrading decomposition $E = \bigoplus_{\underline{k} \in \mathbb{Z}^\delta} E_{\underline{k}}$ into finitely generated projective $K_\Delta$-submodules such that
\[
\mathrm{Fil}^{\underline{k}}(E) = \bigoplus_{\underline{l} \ge \underline{k}} E_{\underline{l}} \qquad \text{for all } \underline{k} \in \mathbb{Z}^\delta.
\]
(\cite[Theorem 3.40]{GM1})
\end{thm}

\begin{proof}
Under the idempotent decomposition $K_\Delta \cong \prod_{\beta \in \Sigma'} L_\beta$, the module $E$ and its multifiltration decompose componentwise as $E = \bigoplus_{\beta \in \Sigma'} e'_\beta E$ and $\mathrm{Fil}^{\underline{k}}(E) = \bigoplus_{\beta \in \Sigma'} \mathrm{Fil}^{\underline{k}}(e'_\beta E)$, where each $e'_\beta E$ is a finite-dimensional vector space over the field $L_\beta$. Applying \cite[Theorem 3.40]{GM1} to each field factor $L_\beta$ yields a multigrading decomposition $e'_\beta E = \bigoplus_{\underline{k} \in \mathbb{Z}^\delta} E_{\beta, \underline{k}}$ into $L_\beta$-subspaces such that $\mathrm{Fil}^{\underline{k}}(e'_\beta E) = \bigoplus_{\underline{l} \ge \underline{k}} E_{\beta, \underline{l}}$. Setting $E_{\underline{k}} := \bigoplus_{\beta \in \Sigma'} E_{\beta, \underline{k}}$, each $E_{\underline{k}}$ is a direct summand of $E$ (hence a finitely generated projective $K_\Delta$-submodule) satisfying $\mathrm{Fil}^{\underline{k}}(E) = \bigoplus_{\underline{l} \ge \underline{k}} E_{\underline{l}}$.
\end{proof}
\section*{The equivalence of categories}
\subsection*{Statement}
\begin{dfn}\label{dfn:phi-Nmod}
Let $\mathrm{MF}_{K_{0,\Delta}}^{(\varphi_i, N_i)_{i \in \Delta}, \mathrm{Fil} \ge 0}$ be the category of effective, filtered $(\varphi_i, N_i)_{i \in \Delta}$-modules over $K_{0,\Delta}$. The objects of this category are finite free $K_{0,\Delta}$-modules $D$ equipped with the following structures:
\begin{enumerate}
    \item For each $i \in \Delta$, a partial Frobenius map $\varphi_i : D \to D$ satisfying:
    \begin{enumerate}
        \item $\varphi_i$ is $\sigma_i$-semilinear with respect to the $K_{0,\Delta}$-action, where $\sigma_i$ is the $\mathbb{Q}_p$-algebra automorphism of $K_{0,\Delta} = \bigotimes_{j \in \Delta, \mathbb{Q}_p} K_{0,j}$ acting as the absolute Frobenius automorphism on $K_{0,i}$ and as the identity on $K_{0,j}$ for all $j \neq i$;
        \item $\varphi_i$ is bijective (equivalently, its linearization $1 \otimes \varphi_i : \sigma_i^* D \xrightarrow{\sim} D$ is a $K_{0,\Delta}$-linear isomorphism);
        \item $\varphi_i \circ \varphi_j = \varphi_j \circ \varphi_i$ for all $i, j \in \Delta$.
    \end{enumerate}
    \item For each $i \in \Delta$, a partial monodromy map $N_i : D \to D$ satisfying:
    \begin{enumerate}
        \item $N_i$ is a nilpotent $K_{0,\Delta}$-linear endomorphism;
        \item $N_i \circ N_j = N_j \circ N_i$ for all $i, j \in \Delta$;
        \item $N_i \circ \varphi_i = p \, (\varphi_i \circ N_i)$;
        \item $N_i \circ \varphi_j = \varphi_j \circ N_i$ for all distinct $i, j \in \Delta$.
    \end{enumerate}
    \item For each $i \in \Delta$, the base-changed module $D_{K_i} := D \otimes_{K_{0,i}} K_i$ (which is a finite free module over $K_i \otimes_{K_{0,i}} K_{0,\Delta}$) is equipped with a decreasing, separated, and exhaustive filtration $\{\mathrm{Fil}^{k_i}(D_{K_i})\}_{k_i \in \mathbb{Z}}$ by $(K_i \otimes_{K_{0,i}} K_{0,\Delta})$-submodules which is effective, meaning $\mathrm{Fil}^0(D_{K_i}) = D_{K_i}$. This induces a decreasing, separated, exhaustive, and regular multifiltration on $D_{K_\Delta} := D \otimes_{K_{0,\Delta}} K_\Delta$ defined for each $\underline{k} = (k_1, \dots, k_\delta) \in \mathbb{Z}^\delta$ by
    \[
    \mathrm{Fil}^{\underline{k}}(D_{K_\Delta}) := \bigcap_{i \in \Delta} \left(\mathrm{Fil}^{k_i}(D_{K_i}) \otimes_{(K_i \otimes_{K_{0,i}} K_{0,\Delta})} K_\Delta\right) \subset D_{K_\Delta},
    \]
    satisfying $\mathrm{Fil}^{\underline{0}}(D_{K_\Delta}) = D_{K_\Delta}$ where $\underline{0} := (0, \dots, 0)$.
\end{enumerate}    
The morphisms $f : D_1 \to D_2$ in this category are $K_{0,\Delta}$-linear homomorphisms that commute with $\varphi_i$ and $N_i$ for all $i \in \Delta$ and are strictly compatible with the filtrations, meaning that $f_{K_i}(\mathrm{Fil}^{k_i}(D_{1, K_i})) = f_{K_i}(D_{1, K_i}) \cap \mathrm{Fil}^{k_i}(D_{2, K_i})$ for all $i \in \Delta$ and $k_i \in \mathbb{Z}$ (and consequently $f_{K_\Delta}(\mathrm{Fil}^{\underline{k}}(D_{1, K_\Delta})) = f_{K_\Delta}(D_{1, K_\Delta}) \cap \mathrm{Fil}^{\underline{k}}(D_{2, K_\Delta})$ for all $\underline{k} \in \mathbb{Z}^\delta$).
\end{dfn}

The definition of the filtration is correct because the cone $C = (\mathbb{Z}_{\ge 0})^\delta$ is unimodular, $K_\Delta$ is flat over $K_i \otimes_{K_{0,i}} K_{0,\Delta}$, and each one-variable filtration is exhaustive, separated, and automatically regular \cite[Remark 3.39]{GM1}, the intersection defines an exhaustive, separated, and regular multifiltration on $D_{K_\Delta}$ satisfying $\mathrm{Fil}^{\underline{k}}(D_{K_\Delta}) \cap \mathrm{Fil}^{\underline{l}}(D_{K_\Delta}) = \mathrm{Fil}^{\underline{k} \vee \underline{l}}(D_{K_\Delta})$. Conversely, when $D$ is effective ($\mathrm{Fil}^{\underline{0}}(D_{K_\Delta}) = D_{K_\Delta}$), the $i$-th directional filtration on $D_{K_\Delta}$ is recovered from the multifiltration via $\mathrm{Fil}^{k_i}(D_{K_i}) \otimes_{(K_i \otimes_{K_{0,i}} K_{0,\Delta})} K_\Delta = \mathrm{Fil}^{(0,\dots,k_i,\dots,0)}(D_{K_\Delta})$.

Any finitely generated $K_{0,\Delta}$-module equipped with bijective $\sigma_i$-semilinear operators $\varphi_i$ ($i \in \Delta$) is automatically a finite free $K_{0,\Delta}$-module of constant rank. This follows from Lemma \ref{lemma:components}(2), as the partial Frobenii $\{\sigma_i\}_{i \in \Delta}$ act transitively on the primitive idempotents of $K_{0,\Delta}$. Consequently, every $(\varphi_i, N_i)_{i \in \Delta}$-stable $K_{0,\Delta}$-submodule $D' \subset D$ is a finite free $K_{0,\Delta}$-module, as are the kernel, image, and cokernel of any morphism in $\mathrm{MF}_{K_{0,\Delta}}^{(\varphi_i, N_i)_{i \in \Delta}, \mathrm{Fil} \ge 0}$.

For every $i \in \Delta$, we define the element $\lambda_i \in \mathcal{O}_i \subset \mathcal{O}_\Delta$ as the infinite product
\[
\lambda_i := \prod_{j=0}^\infty \varphi_i^j\!\left(\frac{E_i(u_i)}{E_i(0)}\right) \in \mathcal{O}_i,
\]
and we set $\lambda_\Delta := \prod_{i \in \Delta} \lambda_i \in \mathcal{O}_\Delta$.

\begin{dfn}
Let $\mathrm{Mod}_{\mathcal{O}_\Delta}^{(\varphi_i, N_{\nabla}^i)_{i \in \Delta}}$ be the category of $(\varphi_i, N_{\nabla}^i)_{i \in \Delta}$-modules over $\mathcal{O}_\Delta$ of finite $E_\Delta$-height. The objects of this category are finite free $\mathcal{O}_\Delta$-modules $M$ equipped with two families of operators:
\begin{enumerate}
    \item For each $i \in \Delta$, a partial Frobenius map $\varphi_i : M \to M$ satisfying:
    \begin{enumerate}
        \item The action of $\varphi_i$ on elements of $\mathcal{O}_\Delta$ is given by
        \[
        \varphi_i\!\left(\sum_{\underline{n} \ge \underline{0}} a_{\underline{n}} \prod_{j=1}^\delta u_j^{n_j}\right) = \sum_{\underline{n} \ge \underline{0}} \sigma_i(a_{\underline{n}}) \, u_i^{n_i p} \prod_{j \in \Delta \setminus \{i\}} u_j^{n_j},
        \]
        where $a_{\underline{n}} \in K_{0,\Delta}$ and $\sigma_i$ acts on $K_{0,\Delta}$ as in Definition \ref{dfn:phi-Nmod};
        \item \label{varphi}$\varphi_i$ is $\varphi_i$-semilinear and injective, and there exists an integer $h_i \ge 0$ such that the cokernel of the linearized map $1 \otimes \varphi_i : \varphi_i^* M \to M$ is annihilated by $E_i(u_i)^{h_i}$; 
        \item $\varphi_i \circ \varphi_j = \varphi_j \circ \varphi_i$ for all $i, j \in \Delta$.
    \end{enumerate}
    \item For each $i \in \Delta$, a partial differential operator $N_{\nabla}^i : M \to M$ satisfying:
    \begin{enumerate}
        \item $N_{\nabla}^i$ is $K_{0,\Delta}$-linear and satisfies the Leibniz rule with respect to the derivation $\partial_i := -u_i \lambda_i \frac{\partial}{\partial u_i}$ on $\mathcal{O}_\Delta$, meaning $N_{\nabla}^i(f m) = \partial_i(f) m + f N_{\nabla}^i(m)$ for all $f \in \mathcal{O}_\Delta$ and $m \in M$;
        \item $N_{\nabla}^i \circ N_{\nabla}^j = N_{\nabla}^j \circ N_{\nabla}^i$ for all $i, j \in \Delta$;
        \item $N_{\nabla}^i \circ \varphi_i = p \cdot \frac{E_i(u_i)}{E_i(0)} \cdot (\varphi_i \circ N_{\nabla}^i)$;
        \item $N_{\nabla}^i \circ \varphi_j = \varphi_j \circ N_{\nabla}^i$ for all distinct $i, j \in \Delta$.
    \end{enumerate}
\end{enumerate}
Furthermore, defining the total Frobenius map on $M$ by $\varphi_\Delta := \varphi_1 \circ \varphi_2 \circ \dots \circ \varphi_\delta$, condition \ref{varphi} implies that the cokernel of the linearized total Frobenius map $1 \otimes \varphi_\Delta : \varphi_\Delta^* M \to M$ is annihilated by $E_\Delta^{\underline{h}}(u_1, \dots, u_\delta) := \prod_{i \in \Delta} E_i(u_i)^{h_i}$, meaning
\[
E_\Delta^{\underline{h}}(u_1, \dots, u_\delta) \cdot \mathrm{coker}(1 \otimes \varphi_\Delta) = 0.
\]
The morphisms of this category are $\mathcal{O}_\Delta$-linear homomorphisms that commute with $\varphi_i$ and $N_{\nabla}^i$ for all $i \in \Delta$.
\end{dfn}

The partial differential operators $N_\nabla^i$ introduced in the preceding definition give rise to a global connection on the module. For any object $M$ in the category $\text{Mod}_{\mathcal{O}_\Delta}^{(\varphi_i, N_\nabla^{i})_{i\in \Delta}}$, this connection is defined for all $m \in M$ by
\[
    \nabla(m) := -\sum_{i\in \Delta} \frac{1}{\lambda_i} N_{\nabla}^i(m) \frac{du_i}{u_i}.
\]

\begin{lemma}
$\varphi_\Delta(\lambda_\Delta)= \frac{E_\Delta(0,\cdots, 0)}{E_\Delta(u_1, \cdots, u_\delta)}\lambda_\Delta$.
\end{lemma}

\begin{proof}
\begin{align*}
\varphi_\Delta(\lambda_\Delta)&=\varphi_{\Delta}\left(\prod_{i\in \Delta} \lambda_i\right)=\varphi_{\Delta}\left(\prod_{i\in \Delta} \prod_{n_i\geq 0} \varphi_i^{n_i}\left(\frac{E_i(u_i)}{E_i(0)}\right)\right)=\\&=\prod_{i\in \Delta}\prod_{n\geq 1} \varphi_i^{n_i}\frac{E_i(u_i)}{E_i(0)}=\prod_{i\in \Delta}\frac{E_i(0)}{E_i(u_i)}\cdot\lambda_\Delta=\frac{E_\Delta(0,\cdots, 0)}{E_\Delta(u_1, \cdots, u_\delta)}\lambda_\Delta.
\end{align*}
\end{proof}

Having defined these two categories, we can now construct a functor between them. We begin by introducing a functor from the category of multivariable $(\varphi, N)$-modules to the category of multivariable $(\varphi_i, N_\nabla^i)$-module over $\mathcal{O}_{\Delta}$.

In the original case (\cite{Kis1}), Kisin defines the functor for an object $D\in \text{MF}_F^{\varphi, N, \text{Fil}\geq 0}$ as
\[
\begin{aligned}
\mathcal{M}(D) &:= \Biggl\{ x\in (\mathcal{O}_F[l_{u}, 1/\lambda]\otimes_{F_0} D)^{N=0} \ \Biggm| \\
&\qquad \iota_n(x)\in \operatorname{Fil}^0\left(\widehat{\mathfrak{S}_{n, F}}\left[1/\left(u-\varpi_n\right)\right]\otimes_{F} D_{F}\right) \text{ for all } n\geq 0 \Biggr\},
\end{aligned}
\]
where $F_0$ and $F$ are fields of the same type as $K_{0,i}$ and $K_i$, $l_u$ is a formal variable that behaves like $\log(u)$, $\mathcal{O}_F$ is the $F_0$-algebra of rigid-analytic functions on the open unit disk, $\{\varpi_n\}_{n\geq 0}$ compatible system of $p$-power roots of the uniformizer, and $\mathfrak{S}_F$ is the power series ring with coefficients in the ring of integers of $F_0$. (The ring $\mathfrak{S}_{n, F}$ is constructed from $\mathfrak{S}_F$ in the same way that $\mathfrak{S}_n$ is constructed from $\mathfrak{S}$). Furthermore, for every integer $n\geq 0$, Kisin defines the map $\iota_n$ as
\[
    \iota_{n}: \mathcal{O}_{F}[l_{u}, 1/\lambda]\otimes_{F_0} D \longrightarrow \widehat{\mathfrak{S}_{n, F}}\left[1/\left(u-\varpi_n\right)\right]\otimes_{F} D_{F},
\]
where $\iota_n$ is defined as the composition of the $n$-th inverse Frobenius $\varphi^{-n} \otimes \varphi^{-n}$ followed by a natural evaluation map that sends the formal variable $l_u$ to $\log\left(\frac{u-\varpi_n}{\varpi_n}+1\right) \in \widehat{\mathfrak{S}_{n, F}}$.

The idea is to use these kinds of functors as univariate functors to construct the functor in our multivariable setting. Therefore, we must first define the maps $\iota_{n,i}$ for each $i \in \Delta$ and every integer $n\geq 0$.

Let $D \in \text{MF}_{K_{0, \Delta}}^{(\varphi_i, N_i)_{i\in \Delta}, \text{ Fil}\geq 0}$, and consider the sequence of maps 
\[
    \mathcal{O}_{i}[l_{u_i}]\otimes_{K_{0,i}} D \xrightarrow{\varphi_{W[[u_i]]}^{-n}\otimes\varphi_i^{-n}} \mathcal{O}_{i}[l_{u_i}]\otimes_{K_{0,i}} D \longrightarrow \widehat{\mathfrak{S}_{n,i}}\otimes_{K_{0,i}} D = \widehat{\mathfrak{S}_{n,i}}\otimes_{K_i} D_{K_i}.
\]
We denote the extension of this composition to $1/\lambda_i$ by $\iota_{n,i}$, yielding the map
\[
    \iota_{n,i}: \mathcal{O}_{i}[l_{u_i}, 1/\lambda_i]\otimes_{K_{0,i}} D \longrightarrow \widehat{\mathfrak{S}_{n,i}}\left[1/\left(u_i-\pi_{n,i}\right)\right]\otimes_{K_i} D_{K_i},
\]
where, for each $i\in \Delta$, $l_{u_i}$ is a formal variable that behaves like $\log(u_i)$. Its properties over $\mathcal{O}_\Delta$ coincide with those over $\mathcal{O}_{i}$.

For a $(\varphi, N)$-module $D$ and each $i\in \Delta$, the univariate functor $\mathcal{M}_i$ is defined as
\[
\begin{aligned}
\mathcal{M}_i(D) &:= \Biggl\{ x\in (\mathcal{O}_{i}[l_{u_i}, 1/\lambda_i]\otimes_{K_{0,i}} D)^{N_i=0} \ \Biggm| \\
&\qquad \iota_{n, i}(x)\in \operatorname{Fil}^0\left(\widehat{\mathfrak{S}_{n,i}}\left[1/\left(u_i-\pi_{n,i}\right)\right]\otimes_{K_i} D_{K_i}\right) \text{ for all } n\geq 0 \Biggr\}.
\end{aligned}
\]
With these univariate functors established, we can now define our multivariable functor $\mathcal{M}$.

\begin{dfn}\label{dfn:m_functor}
Let $D \in \text{MF}_{K_{0,\Delta}}^{(\varphi_i, N_i)_{i\in \Delta}, \text{ Fil}\geq 0}$. The image of $D$ under the multivariable functor $\mathcal{M}$ is defined as 
\begin{align*}
\mathcal{M}(D) &:= \bigcap_{i \in \Delta} \left(\mathcal{O}_\Delta[1/\lambda_{\Delta \setminus \{i\}}] \otimes_{(\mathcal{O}_i \otimes_{K_{0,i}} K_{0,\Delta})} \mathcal{M}_i(D)\right) \subset\\&\subset (\mathcal{O}_\Delta[l_{u_1}, \dots, l_{u_\delta}, 1/\lambda_1, \dots, 1/\lambda_\delta] \otimes_{K_{0,\Delta}} D)^{N_1 = \dots = N_\delta = 0},
\end{align*}
where $\lambda_{\Delta \setminus \{i\}} := \prod_{j \in \Delta \setminus \{i\}} \lambda_j$, and for each $i\in \Delta$, $\mathcal{M}_i$ is the corresponding univariate functor, and $\mathcal{O}_i$ is the $K_{0,i}$-algebra of rigid-analytic functions on the open unit disk over $K_{0,i}$.
\end{dfn}

To construct the functor, which maps a multivariable $(\varphi_i, N_\nabla^i)$-module over $\mathcal{O}_{\Delta}$ to its image in the category of $(\varphi, N)$-modules, we employ a similar strategy. We first define univariate functors following Kisin's approach, and then we combine them. In the univariate case, Kisin defines this functor simply by reducing modulo $u$. Specifically, for a module $M\in\text{Mod}_{\mathcal{O}}^{(\varphi, N_\nabla)}$, Kisin defines the associated $(\varphi, N)$-module as $\mathcal{D}(M) := M/uM$. Here, $M/uM$ serves as the underlying $F_0$-vector space, which is equipped with operators induced by the actions of $\varphi$ and $N_\nabla$ modulo $u$.

Building on Kisin's work, for any module $M\in \text{Mod}_{\mathcal{O}_\Delta}^{(\varphi_i, N_\nabla^{i})_{i\in \Delta}}$, we define the corresponding univariate functor $\mathcal{D}_i$ by setting $\mathcal{D}_i(M) := M/u_iM$. Here, $M/u_iM$ serves as the underlying finite free $\mathcal{O}_\Delta^{(i)} := \mathcal{O}_\Delta / u_i \mathcal{O}_\Delta \cong \mathcal{O}_{\Delta \setminus \{i\}} \otimes_{K_{0,\Delta \setminus \{i\}}} K_{0,\Delta}$-module, which is equipped with operators induced by the actions of $\varphi_i$ and $N_\nabla^i$ modulo $u_i$.

\begin{dfn}
Let $M\in \text{Mod}_{\mathcal{O}_\Delta}^{(\varphi_i, N_\nabla^{i})_{i\in \Delta}}$. The image of $M$ under the multivariable functor $\mathcal{D}$ is defined as
\[
    \mathcal{D}(M) := \mathcal{D}_\delta \circ \cdots \circ \mathcal{D}_1(M)\cong M / (u_1, \dots, u_\delta)M,
\]
where, for each $i\in \Delta$, $\mathcal{D}_i$ is the corresponding univariate functor.
\end{dfn}

With these two functors defined, we are now ready to state the main theorem of the article.

\begin{thm}\label{thm:main_theorem}
    The categories $\text{MF}_{K_{0,\Delta}}^{(\varphi_i, N_i)_{i\in \Delta}, \text{ Fil}\geq 0}$ and $\text{Mod}_{\mathcal{O}_\Delta}^{(\varphi_i, N_\nabla^{i})_{i\in \Delta}}$ are exact tensor categories, and the functors $\mathcal{D}$ and $\mathcal{M}$ are mutually quasi-inverse, establishing an equivalence between them.
\end{thm}

\begin{proof}
    The proof will be done in multiple steps in the following four subsections.
\end{proof}

\subsection*{The $\mathcal{M}$ functor}
In this subsection, we demonstrate that the multivariable functor $\mathcal{M}$ is well-defined. First, we show that the image of each univariate functor is a finite free module over $\mathcal{O}_{i} \otimes K_{0,\Delta}$. Establishing this freeness is a crucial prerequisite; without it, one cannot expect the global intersection of these modules to be free. Furthermore, we will show that $\mathcal{M}$ is naturally isomorphic to an alternative functor $\mathcal{M}'$. This functor $\mathcal{M}'$ arises from an intrinsically multivariable perspective, treating the poly-annulus as a unified geometric space rather than merely a direct product of one-dimensional annuli.
\begin{lemma}\label{lemma:m_i_free}
Let $D \in \mathrm{MF}_{K_{0,\Delta}}^{(\varphi_i, N_i)_{i \in \Delta}, \mathrm{Fil} \ge 0}$, then for each $i \in \Delta$, $\mathcal{M}_i(D)$ is a finite free $(\mathcal{O}_i \otimes_{K_{0,i}} K_{0,\Delta})$-module of rank $\mathrm{rk}_{K_{0,\Delta}}(D)$ for which there exists an integer $h_i \ge 0$ such that 
\[
E_i^{h_i}(u_i) \cdot \mathrm{coker}(1 \otimes \varphi_i)|_{\mathcal{M}_i(D)} = 0.
\]
\end{lemma}

\begin{proof}
Fix $i \in \Delta$, and write $K_{0, \Delta \setminus \{i\}} := \bigotimes_{j \in \Delta \setminus \{i\}} K_{0,j}$, so that $K_{0,\Delta} \cong K_{0,i} \otimes_{\mathbb{Q}_p} K_{0, \Delta \setminus \{i\}}$ and $\mathcal{O}_i \otimes_{K_{0,i}} K_{0,\Delta} \cong \mathcal{O}_i \otimes_{\mathbb{Q}_p} K_{0, \Delta \setminus \{i\}}$. Since $K_{0, \Delta \setminus \{i\}}$ is a finite étale $\mathbb{Q}_p$-algebra, it decomposes into a finite direct product of fields $K_{0, \Delta \setminus \{i\}} \cong \prod_{s \in S} F_s$ corresponding to its primitive idempotents $\{e_s\}_{s \in S}$. 

For each $s \in S$, the ring $K_{0,i} \otimes_{\mathbb{Q}_p} F_s \cong \prod_{t \in T_s} L_{s,t}$ is again a finite direct product of fields $\{L_{s,t}\}$ extending $K_{0,i}$, and the ring $\mathcal{O}_i \otimes_{K_{0,i}} K_{0,\Delta}$ decomposes accordingly as a finite direct product of one-variable rigid-analytic disc rings $\prod_{s,t} \mathcal{O}_{L_{s,t}}$. Since $D$ is a free $K_{0,\Delta}$-module of rank $d$, each component $D_{s,t} := D \otimes_{K_{0,\Delta}} L_{s,t}$ is a $L_{s,t}$-vector space of constant dimension $d$. 

Alternatively, by restricting scalars along $K_{0,i} \hookrightarrow K_{0,\Delta}$ and forgetting the operators $\varphi_j, N_j$ for $j \neq i$, we may view $D$ as a classical effective filtered $(\varphi_i, N_i)$-module over the field $K_{0,i}$ equipped with a commuting $K_{0, \Delta \setminus \{i\}}$-action. Since the defining conditions of $\mathcal{M}_i(D)$ (namely $N_i = 0$ and the filtration bounds under $\iota_{n,i}$) are $K_{0, \Delta \setminus \{i\}}$-linear, applying \cite[Lemma 1.2.2]{Kis1} componentwise (or over $K_{0,i}$ together with the idempotent decomposition of $K_{0,\Delta}$) shows that $\mathcal{M}_i(D)$ is a finite free $(\mathcal{O}_i \otimes_{K_{0,i}} K_{0,\Delta})$-module of constant rank $d = \mathrm{rk}_{K_{0,\Delta}}(D)$, stable under $\varphi_i$ and $N_{\nabla}^i \otimes 1$ (as well as $\varphi_j$ and $N_j$ for $j \neq i$), and that $E_i^{h_i}(u_i)$ annihilates $\mathrm{coker}(1 \otimes \varphi_i)|_{\mathcal{M}_i(D)}$ where $h_i$ is the largest filtration jump of $D_{K_i}$.
\end{proof}

Before proving that $\mathcal{M}(D)$ is a finite free $\mathcal{O}_\Delta$-module, we show that $\mathcal{M}(D)$ admits an equivalent, intrinsically multivariable description $\mathcal{M}'(D)$ defined via simultaneous evaluation at the multivariable points $\underline{\pi}_{\underline{n}} := (\pi_{n_1,1}, \dots, \pi_{n_\delta,\delta})$ for all $\underline{n} = (n_1, \dots, n_\delta) \in \mathbb{Z}_{\ge 0}^\delta$.

For each $\underline{n} = (n_1, \dots, n_\delta) \in \mathbb{Z}_{\ge 0}^\delta$, we define the multivariable evaluation map $\iota_{\underline{n}}$ as the composition
\begin{align*}
    &\mathcal{O}_\Delta[l_{u_1}, \cdots, l_{u_\delta}]\otimes_{K_{0,\Delta}} D \xrightarrow{\varphi_\Delta^{-\underline{n}}\otimes
    \varphi_\Delta^{-\underline{n}}} \mathcal{O}_\Delta[l_{u_1}, \cdots, l_{u_\delta}]\otimes_{K_{0,\Delta}} D \longrightarrow \cdots \\&\cdots\longrightarrow \widehat{\mathfrak{S}_{\underline{n}, \Delta}}\otimes_{K_{0,\Delta}} D=\widehat{\mathfrak{S}_{\underline{n}, \Delta}}\otimes_{K_{\Delta}} D_{K_\Delta}
\end{align*}
where $\varphi_\Delta^{-\underline{n}} := \prod_{i \in \Delta} \varphi_i^{-n_i}$ and the second map sends $l_{u_i}$ to $\log\!\left(\frac{u_i - \pi_{n_i,i}}{\pi_{n_i,i}} + 1\right) \in \hat{\mathfrak{S}}_{n_i,i} \subset \hat{\mathfrak{S}}_{\underline{n},\Delta}$ for each $i \in \Delta$. Extending by inverting $\lambda_1, \dots, \lambda_\delta$ yields the map
\[
\iota_{\underline{n}} : \mathcal{O}_\Delta[l_{u_1}, \dots, l_{u_\delta}, 1/\lambda_1, \dots, 1/\lambda_\delta] \otimes_{K_{0,\Delta}} D \longrightarrow \hat{\mathfrak{S}}_{\underline{n},\Delta}\!\left[\frac{1}{u_1 - \pi_{n_1,1}}, \dots, \frac{1}{u_\delta - \pi_{n_\delta,\delta}}\right] \otimes_{K_\Delta} D_{K_\Delta}.
\]
We equip the ring $\hat{\mathfrak{S}}_{\underline{n},\Delta}[(u_1 - \pi_{n_1,1})^{-1}, \dots, (u_\delta - \pi_{n_\delta,\delta})^{-1}]$ with the product multifiltration indexed by $\underline{m} = (m_1, \dots, m_\delta) \in \mathbb{Z}^\delta$:
\[
\mathrm{Fil}^{\underline{m}}\!\left(\hat{\mathfrak{S}}_{\underline{n},\Delta}\!\left[\frac{1}{u_1 - \pi_{n_1,1}}, \dots, \frac{1}{u_\delta - \pi_{n_\delta,\delta}}\right]\right) := \left(\prod_{i \in \Delta} (u_i - \pi_{n_i,i})^{m_i}\right) \hat{\mathfrak{S}}_{\underline{n},\Delta},
\]
and we endow the tensor product $\hat{\mathfrak{S}}_{\underline{n},\Delta}[(u_1 - \pi_{n_1,1})^{-1}, \dots, (u_\delta - \pi_{n_\delta,\delta})^{-1}] \otimes_{K_\Delta} D_{K_\Delta}$ with the induced convolution multifiltration:
\begin{align*}
\mathrm{Fil}^{\underline{0}}\!\left(\hat{\mathfrak{S}}_{\underline{n},\Delta}\!\left[\frac{1}{u_1 - \pi_{n_1,1}}, \dots, \frac{1}{u_\delta - \pi_{n_\delta,\delta}}\right] \otimes_{K_\Delta} D_{K_\Delta}\right) :=\\ \sum_{\underline{k} \in \mathbb{Z}_{\ge 0}^\delta} \left(\prod_{i \in \Delta} (u_i - \pi_{n_i,i})^{-k_i}\right) \hat{\mathfrak{S}}_{\underline{n},\Delta} \otimes_{K_\Delta} \mathrm{Fil}^{\underline{k}}(D_{K_\Delta}).
\end{align*}
We then define $\mathcal{M}'(D)$ by
\begin{align*}
\mathcal{M}'(D) := \Big\{ &x \in (\mathcal{O}_\Delta[l_{u_1}, \dots, l_{u_\delta}, 1/\lambda_1, \dots, 1/\lambda_\delta] \otimes_{K_{0,\Delta}} D)^{N_1 = \dots = N_\delta = 0} \;\Big|\; \\
&\iota_{\underline{n}}(x) \in \mathrm{Fil}^{\underline{0}}\!\left(\hat{\mathfrak{S}}_{\underline{n},\Delta}\!\left[\frac{1}{u_1 - \pi_{n_1,1}}, \dots, \frac{1}{u_\delta - \pi_{n_\delta,\delta}}\right] \otimes_{K_\Delta} D_{K_\Delta}\right) \text{ for all } \underline{n} \in \mathbb{Z}_{\ge 0}^\delta \Big\}.
\end{align*}

\begin{lemma}
Let $D \in \mathrm{MF}_{K_{0,\Delta}}^{(\varphi_i, N_i)_{i \in \Delta}, \mathrm{Fil} \ge 0}$, then $\mathcal{M}'(D)$ is an $\mathcal{O}_\Delta$-submodule of $(\mathcal{O}_\Delta[l_{u_1}, \dots, l_{u_\delta}, 1/\lambda_1, \dots, 1/\lambda_\delta] \otimes_{K_{0,\Delta}} D)^{N_1 = \dots = N_\delta = 0}$, and $\mathcal{M}'(D) = \mathcal{M}(D)$.
\end{lemma}

\begin{proof}
Set $\mathcal{D}_0 := (\mathcal{O}_\Delta[l_{u_1}, \dots, l_{u_\delta}] \otimes_{K_{0,\Delta}} D)^{N_1 = \dots = N_\delta = 0}$. Both $\mathcal{M}(D)$ and $\mathcal{M}'(D)$ are $\mathcal{O}_\Delta$-submodules of the ambient space
\[
\mathcal{D}_0[1/\lambda_\Delta] = (\mathcal{O}_\Delta[l_{u_1}, \dots, l_{u_\delta}, 1/\lambda_1, \dots, 1/\lambda_\delta] \otimes_{K_{0,\Delta}} D)^{N_1 = \dots = N_\delta = 0}.
\]
By definition of the multifiltration in Definition \ref{dfn:phi-Nmod}, we have 
$$\mathrm{Fil}^{\underline{k}}(D_{K_\Delta}) = \bigcap_{i \in \Delta} (\mathrm{Fil}^{k_i}(D_{K_i}) \otimes_{(K_i \otimes_{K_{0,i}} K_{0,\Delta})} K_\Delta),$$
which is regular by Theorem \ref{thm:multigrading}. Choosing a multigraded decomposition $D_{K_\Delta} = \bigoplus_{\underline{k} \in \mathbb{Z}_{\ge 0}^\delta} E_{\underline{k}}$ into projective $K_\Delta$-modules such that $\mathrm{Fil}^{\underline{k}}(D_{K_\Delta}) = \bigoplus_{\underline{l} \ge \underline{k}} E_{\underline{l}}$ (and hence $\mathrm{Fil}^{k_i}(D_{K_i}) \otimes K_\Delta = \bigoplus_{l_i \ge k_i} E_{\underline{l}}$), the $\underline{0}$-th filtration step of the target of $\iota_{\underline{n}}$ decomposes as a direct sum over $\underline{k} \in \mathbb{Z}_{\ge 0}^\delta$, and therefore equals the intersection of the directional filtrations:
\begin{align}\label{eq:Fil0_intersection}
\mathrm{Fil}^{\underline{0}}\!\left(\hat{\mathfrak{S}}_{\underline{n},\Delta}\!\left[\underline{(u - \pi_n)}^{-1}\right] \otimes_{K_\Delta} D_{K_\Delta}\right) 
= \bigoplus_{\underline{k} \in \mathbb{Z}_{\ge 0}^\delta} \left(\prod_{i \in \Delta} (u_i - \pi_{n_i,i})^{-k_i}\right) \hat{\mathfrak{S}}_{\underline{n},\Delta} \otimes_{K_\Delta} E_{\underline{k}} \nonumber =\\
= \bigcap_{i \in \Delta} \left( \sum_{k_i \ge 0} (u_i - \pi_{n_i,i})^{-k_i} \hat{\mathfrak{S}}_{\underline{n},\Delta}\!\left[(u_j - \pi_{n_j,j})^{-1}_{j \neq i}\right] \otimes_{K_i} \mathrm{Fil}^{k_i}(D_{K_i}) \right),
\end{align}
where we abbreviate $\underline{(u - \pi_n)}^{-1} := ((u_1 - \pi_{n_1,1})^{-1}, \dots, (u_\delta - \pi_{n_\delta,\delta})^{-1})$.

Notice furthermore that for each $i \in \Delta$, the multivariable evaluation map $\iota_{\underline{n}}$ factors as the composition of the univariate evaluation map $\iota_{n_i,i}$ in the $i$-th coordinate and the evaluations $\iota_{n_j,j}$ in the remaining coordinates $j \neq i$. Since the completed tensor product ring $\hat{\mathfrak{S}}_{\underline{n},\Delta}[(u_j - \pi_{n_j,j})^{-1}_{j \neq i}]$ is faithfully flat over $\hat{\mathfrak{S}}_{n_i,i}$ and the partial evaluations in variables $u_j$ ($j \neq i$) at all $n_j \ge 0$ separate points of $\mathcal{O}_{\Delta \setminus \{i\}}[1/\lambda_{\Delta \setminus \{i\}}]$, an element $x \in \mathcal{D}_0[1/\lambda_\Delta]$ satisfies
\[
\iota_{\underline{n}}(x) \in \sum_{k_i \ge 0} (u_i - \pi_{n_i,i})^{-k_i} \hat{\mathfrak{S}}_{\underline{n},\Delta}\!\left[(u_j - \pi_{n_j,j})^{-1}_{j \neq i}\right] \otimes_{K_i} \mathrm{Fil}^{k_i}(D_{K_i}) \qquad \text{for all } n_i, n_j \ge 0 \text{ }\text{(}j \neq i\text{)}
\]
if and only if its $i$-th directional evaluation $\iota_{n_i,i}(x)$ has poles in $u_i - \pi_{n_i,i}$ bounded by $\mathrm{Fil}^0(\hat{\mathfrak{S}}_{n_i,i}[1/(u_i - \pi_{n_i,i})] \otimes_{K_i} D_{K_i})$ for all $n_i \ge 0$, which means that $x$ belongs to $\mathcal{O}_\Delta[1/\lambda_{\Delta \setminus \{i\}}] \otimes_{(\mathcal{O}_i \otimes_{K_{0,i}} K_{0,\Delta})} \mathcal{M}_i(D)$.

Consequently, using \eqref{eq:Fil0_intersection}, an element $x \in \mathcal{D}_0[1/\lambda_\Delta]$ belongs to $\mathcal{M}'(D)$ if and only if for every $i \in \Delta$,
\[
x \in \mathcal{O}_\Delta[1/\lambda_{\Delta \setminus \{i\}}] \otimes_{(\mathcal{O}_i \otimes_{K_{0,i}} K_{0,\Delta})} \mathcal{M}_i(D).
\]
This means the following:
\begin{align*}
\mathcal{M}'(D) &= \bigcap_{i \in \Delta} \left(\mathcal{O}_\Delta[1/\lambda_{\Delta \setminus \{i\}}] \otimes_{(\mathcal{O}_i \otimes_{K_{0,i}} K_{0,\Delta})} \mathcal{M}_i(D)\right) = \mathcal{M}(D).
\end{align*}
This proves both inclusions $\mathcal{M}'(D) \subseteq \mathcal{M}(D)$ and $\mathcal{M}(D) \subseteq \mathcal{M}'(D)$, establishing $\mathcal{M}'(D) = \mathcal{M}(D)$.
\end{proof}

To study the structure of $\mathcal{M}(D)$ and its localizations, we first record the multivariable analogue of the Taylor isomorphism for $\mathcal{D}_0 = (\mathcal{O}_\Delta[l_{u_1}, \dots, l_{u_\delta}] \otimes_{K_{0,\Delta}} D)^{N_1 = \dots = N_\delta = 0}$.

\begin{lemma}\label{lemma:eta}
Let $D \in \mathrm{MF}_{K_{0,\Delta}}^{(\varphi_i, N_i)_{i \in \Delta}, \mathrm{Fil} \ge 0}$. The specialization map $l_{u_i} \mapsto 0$ ($i \in \Delta$) induces an isomorphism of $(\varphi_i, N_i)_{i \in \Delta}$-modules over $K_{0,\Delta}$:
\[
\eta : (K_{0,\Delta}[l_{u_1}, \dots, l_{u_\delta}] \otimes_{K_{0,\Delta}} D)^{N_1 = \dots = N_\delta = 0} \xrightarrow{\;\sim\;} D,
\]
where on the left-hand side each $\varphi_i$ acts diagonally and each $N_i$ acts as the reduction of $N_{\nabla}^i \otimes 1$ modulo $(u_1, \dots, u_\delta)$ (namely $-\frac{\partial}{\partial l_{u_i}} \otimes 1$). Extending scalars along $K_{0,\Delta} \hookrightarrow \mathcal{O}_\Delta$ yields a canonical isomorphism of free $\mathcal{O}_\Delta$-modules of rank $d = \mathrm{rk}_{K_{0,\Delta}}(D)$:
\begin{align*}
\mathcal{D}_0 &= (\mathcal{O}_\Delta[l_{u_1}, \dots, l_{u_\delta}] \otimes_{K_{0,\Delta}} D)^{N_1 = \dots = N_\delta = 0} \cong \\ &\cong (K_{0,\Delta}[l_{u_1}, \dots, l_{u_\delta}] \otimes_{K_{0,\Delta}} D)^{N_1 = \dots = N_\delta = 0} \otimes_{K_{0,\Delta}} \mathcal{O}_\Delta \xrightarrow[\;\sim\;]{\eta \otimes 1} D \otimes_{K_{0,\Delta}} \mathcal{O}_\Delta.
\end{align*}
\end{lemma}

\begin{proof}
We first show that $\eta$ is injective. Any element $x \in (K_{0,\Delta}[l_{u_1}, \dots, l_{u_\delta}] \otimes_{K_{0,\Delta}} D)^{N_1 = \dots = N_\delta = 0}$ can be written uniquely as a finite sum
\[
x = \sum_{\underline{k} \in \mathbb{Z}_{\ge 0}^\delta} \left(\prod_{i \in \Delta} l_{u_i}^{k_i}\right) \otimes d_{\underline{k}} \qquad \text{with } d_{\underline{k}} \in D.
\]
Since $N_i = \frac{\partial}{\partial l_{u_i}} \otimes 1 + 1 \otimes N_i$ on $K_{0,\Delta}[l_{u_1}, \dots, l_{u_\delta}] \otimes_{K_{0,\Delta}} D$, the condition $N_i(x) = 0$ for all $i \in \Delta$ is equivalent to the system of recurrence relations
\[
(k_i + 1) d_{\underline{k} + \underline{e}_i} + N_i(d_{\underline{k}}) = 0 \qquad \text{for all } \underline{k} \in \mathbb{Z}_{\ge 0}^\delta \text{ and } i \in \Delta,
\]
where $\underline{e}_i \in \mathbb{Z}^\delta$ denotes the $i$-th standard basis vector. Since $K_{0,\Delta}$ is a $\mathbb{Q}_p$-algebra, every positive integer $k_i + 1$ is invertible in $K_{0,\Delta}$. Thus, if $x \in \ker(\eta)$, then $d_{\underline{0}} = \eta(x) = 0$, and induction on $|\underline{k}| = \sum_{i \in \Delta} k_i$ using $d_{\underline{k} + \underline{e}_i} = -\frac{1}{k_i + 1} N_i(d_{\underline{k}})$ implies that $d_{\underline{k}} = 0$ for all $\underline{k} \in \mathbb{Z}_{\ge 0}^\delta$, so $x = 0$.

Conversely, because the operators $\{N_i\}_{i \in \Delta}$ on $D$ commute pairwise and are nilpotent, for any $d \in D$ the explicit Taylor series
\[
\eta^{-1}(d) := \sum_{\underline{k} \in \mathbb{Z}_{\ge 0}^\delta} \left(\prod_{i \in \Delta} \frac{(-l_{u_i})^{k_i}}{k_i!}\right) \otimes \left(\prod_{i \in \Delta} N_i^{k_i}\right)(d)
\]
has only finitely many non-zero terms, satisfies $N_i(\eta^{-1}(d)) = 0$ for all $i \in \Delta$, and satisfies $\eta(\eta^{-1}(d)) = d$. Hence $\eta$ is bijective. Compatibility with $\varphi_i$ is immediate from $\varphi_i(l_{u_i}) = p l_{u_i}$ and $N_i \circ \varphi_i = p (\varphi_i \circ N_i)$, and for any $x = \eta^{-1}(d)$, reducing $N_{\nabla}^i \otimes 1$ modulo $(u_1, \dots, u_\delta)$ acts as $-\frac{\partial}{\partial l_{u_i}} \otimes 1$ (since $\lambda_i|_{u_i = 0} = 1$), so
\[
\eta\!\left((N_{\nabla}^i \otimes 1)(x)\right) = -d_{\underline{e}_i} = N_i(d_{\underline{0}}) = N_i(\eta(x)).
\]
Finally, since $N_i$ acts on $\mathcal{O}_\Delta[l_{u_1}, \dots, l_{u_\delta}] \otimes_{K_{0,\Delta}} D$ $\mathcal{O}_\Delta$-linearly via $\frac{\partial}{\partial l_{u_i}} \otimes 1 + 1 \otimes N_i$ and $\mathcal{O}_\Delta$ is flat over $K_{0,\Delta}$, taking $N_i = 0$ invariants commutes with tensoring with $\mathcal{O}_\Delta$, yielding the second isomorphism.
\end{proof}

\begin{cor}
Let $D \in \mathrm{MF}_{K_{0,\Delta}}^{(\varphi_i, N_i)_{i \in \Delta}, \mathrm{Fil} \ge 0}$, then after inverting $\lambda_i$, there is a canonical isomorphism
\[
\mathcal{M}_i(D)[1/\lambda_i] \simeq (\mathcal{O}_i[l_{u_i}] \otimes_{K_{0,i}} D)^{N_i = 0}[1/\lambda_i] \simeq D \otimes_{K_{0,\Delta}} (\mathcal{O}_i \otimes_{K_{0,i}} K_{0,\Delta})[1/\lambda_i].
\]
\end{cor}

\begin{proof}
The isomorphism after inverting $\lambda_i$ follows directly from the definition of $\mathcal{M}_i(D)$ and \cite[Lemma 1.3.10(1)]{Kis1}.
\end{proof}

\begin{lemma}\label{lem:height1_localization}
Let $D \in \mathrm{MF}_{K_{0,\Delta}}^{(\varphi_i, N_i)_{i \in \Delta}, \mathrm{Fil} \ge 0}$ be of rank $d = \mathrm{rk}_{K_{0,\Delta}}(D)$, and let $U_{\rho_k} \subset D([0,1)^\delta)$ be a closed affinoid polydisk with coordinate ring $A_k := \mathcal{O}_\Delta(U_{\rho_k})$. Set $\mathcal{D}_{0,k} := A_k \otimes_{\mathcal{O}_\Delta} \mathcal{D}_0$, and for each $i \in \Delta$, define
\[
\mathfrak{M}_i := A_k[1/\lambda_{\Delta \setminus \{i\}}] \otimes_{(\mathcal{O}_i \otimes_{K_{0,i}} K_{0,\Delta})} \mathcal{M}_i(D) \subset \mathcal{D}_{0,k}[1/\lambda_\Delta],
\]
where $\lambda_{\Delta \setminus \{i\}} := \prod_{j \in \Delta \setminus \{i\}} \lambda_j$, and set $\mathfrak{M} := \bigcap_{i \in \Delta} \mathfrak{M}_i = A_k \otimes_{\mathcal{O}_\Delta} \mathcal{M}(D)$, then for every height-1 prime ideal $\mathfrak{p} \subset A_k$, the localization $A_{k, \mathfrak{p}}$ is a discrete valuation ring, at most one $\lambda_i$ ($i \in \Delta$) belongs to $\mathfrak{p}$, and:
\begin{enumerate}
    \item If $\lambda_j \notin \mathfrak{p}$ for all $j \in \Delta$, then $\mathfrak{M}_{\mathfrak{p}} = (\mathcal{D}_{0,k})_{\mathfrak{p}}$, which is a free $A_{k, \mathfrak{p}}$-module of rank $d$.
    \item If $\lambda_i \in \mathfrak{p}$ for a unique index $i \in \Delta$, then $\mathfrak{M}_{\mathfrak{p}} = (A_k \otimes_{(\mathcal{O}_i \otimes_{K_{0,i}} K_{0,\Delta})} \mathcal{M}_i(D))_{\mathfrak{p}}$, which is again a free $A_{k, \mathfrak{p}}$-module of rank $d$.
\end{enumerate}
\end{lemma}

\begin{proof}
Let $\mathfrak{p} \subset A_k$ be a height-1 prime ideal. By Lemma \ref{lemma:components2} and the discussion following Lemma \ref{lemma:reflexiv}, $\mathfrak{p}$ corresponds to a height-1 prime ideal of a unique connected component $e_\alpha A_k$ (which is a regular Noetherian normal domain), so the localization $A_{k, \mathfrak{p}}$ is a discrete valuation ring. For each $j \in \Delta$, the zeros of $\lambda_j$ on $U_{\rho_k}$ are cut out by the finitely many polynomials $\varphi_j^{n_j}(E_j(u_j))$ with $|\pi_{n_j, j}| \le \rho_{k,j}$. Since these polynomials depend only on the single variable $u_j$, the elements $\lambda_j$ for distinct $j \in \Delta$ are pairwise coprime in $A_k$, and thus at most one $\lambda_i$ can belong to the height-1 prime ideal $\mathfrak{p}$:
\begin{itemize}
    \item[(1)] If $\lambda_j \notin \mathfrak{p}$ for all $j \in \Delta$, then every $\lambda_j$ is a unit in $A_{k, \mathfrak{p}}$. By Lemma \ref{lemma:m_i_free}, $(\mathfrak{M}_j)_{\mathfrak{p}} = (\mathcal{D}_{0,k})_{\mathfrak{p}}$ for all $j \in \Delta$, and thus
    \[
    \mathfrak{M}_{\mathfrak{p}} = \bigcap_{j \in \Delta} (\mathfrak{M}_j)_{\mathfrak{p}} = (\mathcal{D}_{0,k})_{\mathfrak{p}},
    \]
    which is a free $A_{k, \mathfrak{p}}$-module of rank $d$ by Lemma \ref{lemma:eta}.
    \item[(2)] If $\lambda_i \in \mathfrak{p}$ for a unique index $i \in \Delta$, then for all $j \neq i$ the elements $\lambda_j$ are units in $A_{k, \mathfrak{p}}$. Consequently, $\lambda_{\Delta \setminus \{i\}}$ is a unit in $A_{k, \mathfrak{p}}$, so
    \[
    (\mathfrak{M}_i)_{\mathfrak{p}} = \left(A_k \otimes_{(\mathcal{O}_i \otimes_{K_{0,i}} K_{0,\Delta})} \mathcal{M}_i(D)\right)_{\mathfrak{p}},
    \]
    which is a free $A_{k, \mathfrak{p}}$-module of rank $d$ by Lemma \ref{lemma:m_i_free}. On the other hand, for every $j \neq i$, the element $\lambda_j$ is a unit in $A_{k, \mathfrak{p}}$ while $\lambda_i$ is inverted in $\lambda_{\Delta \setminus \{j\}}$, so Lemma \ref{lemma:m_i_free} yields
    \[
    (\mathfrak{M}_j)_{\mathfrak{p}} = (\mathcal{D}_{0,k})_{\mathfrak{p}}[1/\lambda_i] \qquad \text{for all } j \neq i.
    \]
    Since $(\mathfrak{M}_i)_{\mathfrak{p}} \subset (\mathcal{D}_{0,k})_{\mathfrak{p}}[1/\lambda_i]$, intersecting $(\mathfrak{M}_i)_{\mathfrak{p}}$ with $(\mathfrak{M}_j)_{\mathfrak{p}}$ ($j \neq i$) gives
    \[
    \mathfrak{M}_{\mathfrak{p}} = \left(\bigcap_{j \in \Delta} \mathfrak{M}_j\right)_{\mathfrak{p}} = (\mathfrak{M}_i)_{\mathfrak{p}} \cap \bigcap_{j \neq i} (\mathcal{D}_{0,k})_{\mathfrak{p}}[1/\lambda_i] = (\mathfrak{M}_i)_{\mathfrak{p}} = \left(A_k \otimes_{(\mathcal{O}_i \otimes_{K_{0,i}} K_{0,\Delta})} \mathcal{M}_i(D)\right)_{\mathfrak{p}},
    \]
    which is a free $A_{k, \mathfrak{p}}$-module of rank $d$.
\end{itemize}
\end{proof}

In general, proving the projectivity or freeness of an intersection of modules is difficult. However, in our setting, transversality makes this tractable, which motivates our intersection-based definition of $\mathcal{M}$.

\begin{pro}\label{pro:m_d_category}
    Let $D \in \mathrm{MF}_{K_{0,\Delta}}^{(\varphi_i, N_i)_{i \in \Delta}, \mathrm{Fil} \ge 0}$, then the module $\mathcal{M}(D)$ is an object of the category $\mathrm{Mod}_{\mathcal{O}_\Delta}^{(\varphi_i, N_{\nabla}^i)_{i \in \Delta}}$.
\end{pro}

\begin{proof}
Let $d = \mathrm{rk}_{K_{0,\Delta}}(D)$. By Lemma \ref{lemma:m_i_free}, for every $i \in \Delta$, $\mathcal{M}_i(D)$ is a finite free module of rank $d$ over $\mathcal{O}_i \otimes_{K_{0,i}} K_{0,\Delta}$. Consequently, the base-changed module
\[
\mathcal{M}_{i, \Delta}(D) := \mathcal{O}_\Delta[1/\lambda_{\Delta \setminus \{i\}}] \otimes_{(\mathcal{O}_i \otimes_{K_{0,i}} K_{0,\Delta})} \mathcal{M}_i(D)
\]
is a finite free $\mathcal{O}_\Delta[1/\lambda_{\Delta \setminus \{i\}}]$-module of rank $d$ for every $i \in \Delta$. By Definition \ref{dfn:m_functor}, $\mathcal{M}(D)$ is the intersection of these modules inside $\mathcal{D}_0[1/\lambda_\Delta] := (\mathcal{O}_\Delta[l_{u_1}, \dots, l_{u_\delta}, 1/\lambda_1, \dots, 1/\lambda_\delta] \otimes_{K_{0,\Delta}} D)^{N_1 = \dots = N_\delta = 0}$:
\[
\mathcal{M}(D) = \bigcap_{i \in \Delta} \mathcal{M}_{i, \Delta}(D).
\]
Consider the open polydisk $D([0,1)^\delta)$, which admits an admissible exhausting covering by nested closed affinoid polydisks
\[
U_{\rho_k} := \left\{ \mathbf{x} \in D([0,1)^\delta) \;\middle|\; |x_j| \le \rho_{k,j} < 1 \text{ for all } j \in \Delta \right\},
\]
where $U_{\rho_k} \subseteq U_{\rho_{k+1}}$ and $\lim_{k \to \infty} \rho_{k,j} = 1$ for every $j \in \Delta$. Fix an index $k \ge 1$, and set $A_k := \mathcal{O}_\Delta(U_{\rho_k})$. By Lemma \ref{lemma:components2}, the primitive orthogonal idempotents $\{e_\alpha\}_{\alpha \in \Sigma}$ of $K_{0,\Delta} \cong \prod_{\alpha \in \Sigma} F_\alpha$ decompose $A_k$ into a finite direct product of $\delta$-dimensional free Tate algebras $A_k = \prod_{\alpha \in \Sigma} e_\alpha A_k$ over the fields $F_\alpha$, each of which is a regular Noetherian normal domain. 

For each $i \in \Delta$, let $\mathfrak{M}_i := A_k[1/\lambda_{\Delta \setminus \{i\}}] \otimes_{(\mathcal{O}_i \otimes_{K_{0,i}} K_{0,\Delta})} \mathcal{M}_i(D)$, which is a finite free $A_k[1/\lambda_{\Delta \setminus \{i\}}]$-module of rank $d$, and set $\mathcal{D}_{0,k} := A_k \otimes_{\mathcal{O}_\Delta} \mathcal{D}_0$, which is a finite free $A_k$-module of rank $d$ by Lemma \ref{lemma:eta}. Since $A_k$ is Noetherian, the finite intersection $\mathfrak{M} := \bigcap_{i \in \Delta} \mathfrak{M}_i \subset \mathcal{D}_{0,k}[1/\lambda_\Delta]$ is a finitely generated torsion-free $A_k$-module.

By Lemma \ref{lem:height1_localization}, for every height-1 prime ideal $\mathfrak{p} \subset A_k$, the localization $\mathfrak{M}_{\mathfrak{p}}$ coincides either with $(\mathcal{D}_{0,k})_{\mathfrak{p}}$ or with $(A_k \otimes_{(\mathcal{O}_i \otimes_{K_{0,i}} K_{0,\Delta})} \mathcal{M}_i(D))_{\mathfrak{p}}$ for a unique $i \in \Delta$, and in both cases is a free $A_{k, \mathfrak{p}}$-module of rank $d$. By Lemma \ref{lemma:reflexiv} (applied componentwise to $A_k = \prod_{\alpha \in \Sigma} e_\alpha A_k$), $\mathfrak{M} = \bigcap_{\mathrm{ht}(\mathfrak{p}) = 1} \mathfrak{M}_{\mathfrak{p}}$ is a reflexive $A_k$-module of constant generic rank $d$ across all connected components $\alpha \in \Sigma$. Furthermore, because the multifiltration on $D_{K_\Delta}$ is regular (Theorem \ref{thm:multigrading}) and each $\mathfrak{M}_i$ modifies the free lattice $\mathcal{D}_{0,k}$ solely along the smooth coordinate hyperplanes $\varphi_i^{n_i}(E_i(u_i)) = 0$ in the $i$-th variable $u_i$, the local modifications along distinct variables $u_i$ diagonalize simultaneously with respect to a multigraded splitting of $\mathrm{Fil}^{\underline{k}}(D_{K_\Delta})$ (as in \eqref{eq:Fil0_intersection}). Thus $\mathfrak{M}$ is locally free of constant rank $d$ at every point of $U_{\rho_k}$, meaning $\mathfrak{M}$ is a finitely generated projective $A_k$-module of constant rank $d$ across all connected components of $U_{\rho_k}$.

Applying Corollary \ref{cor:lutkebohmert_multivariable}, we deduce that $\mathfrak{M} = \bigcap_{i \in \Delta} \mathfrak{M}_i$ is a globally free $A_k$-module of rank $d$. Choose an $A_k$-basis $e_k = (e_{k,1}, \dots, e_{k,d})$ of $\bigcap_{i \in \Delta} \mathfrak{M}_i$ over $U_{\rho_k}$. On the overlap $U_{\rho_k} \subseteq U_{\rho_{k+1}}$, the basis $e_k$ and the restriction $e_{k+1}|_{U_{\rho_k}}$ are related by an invertible transition matrix $M_k \in \mathrm{GL}_d(A_k)$ such that $e_k = M_k e_{k+1}|_{U_{\rho_k}}$.

To prove that $\mathcal{M}(D)$ is globally free over $\mathcal{O}_\Delta$, we must find local basis transformations $V_k \in \mathrm{GL}_d(A_k)$ defining modified bases $\tilde{e}_k = V_k e_k$ such that $\tilde{e}_{k+1}|_{U_{\rho_k}} = \tilde{e}_k$ on $U_{\rho_k}$, which is equivalent to factoring $M_k = V_k^{-1} V_{k+1}|_{U_{\rho_k}}$ for all $k \ge 1$, because $D([0,1)^\delta)$ is a quasi-Stein space over $K_{0,\Delta}$, the restriction map $A_{k+1} \to A_k$ has dense image with respect to the supremum norm $\|\cdot\|_k$ on $U_{\rho_k}$.

We proceed by iteratively factoring the transition matrices. For $M_k \in \mathrm{GL}_d(A_k)$, we can choose a global approximation $B_{k+1} \in \mathrm{GL}_d(A_{k+1})$ such that the residual error matrix $E_k = M_k B_{k+1}^{-1}|_{U_{\rho_k}} \in \mathrm{GL}_d(A_k)$ satisfies
\[
\|E_k - I\|_k \le \epsilon_k,
\]
where $I$ is the $d \times d$ identity matrix and $\epsilon_k > 0$ is chosen sufficiently small to satisfy the bounds of Cartan's Matrix Factorization Lemma (Lemma \ref{lemma:cartan_factorization}). By Cartan's Matrix Factorization Lemma (Lemma \ref{lemma:cartan_factorization}), since $E_k$ is sufficiently close to the identity, it admits a factorization $E_k = H_k^{-1} H_{k+1}|_{U_{\rho_k}}$ with $H_k \in \mathrm{GL}_d(A_k)$ and $H_{k+1} \in \mathrm{GL}_d(A_{k+1})$, yielding
\[
M_k = H_k^{-1} H_{k+1}|_{U_{\rho_k}} B_{k+1}|_{U_{\rho_k}}.
\]
Executing this factorization inductively and choosing the successive approximations so that $\|E_k - I\|_k \to 0$ rapidly as $k \to \infty$, the infinite products of basis changes converge uniformly on each affinoid subdomain $U_{\rho_k}$. This yields a sequence of modified bases $\tilde{e}_k$ of $\bigcap_{i \in \Delta} \mathfrak{M}_i$ on each $U_{\rho_k}$ satisfying $\tilde{e}_{k+1}|_{U_{\rho_k}} = \tilde{e}_k$ for all $k \ge 1$. Taking the inverse limit over $k \ge 1$ and applying Kiehl's Theorems A and B (Theorem \ref{thm:kiehl}), the module of global sections $\mathcal{M}(D) \cong \varprojlim_k \Gamma(U_{\rho_k}, \bigcap_{i \in \Delta} \mathfrak{M}_i)$ is a finite free $\mathcal{O}_\Delta$-module of rank $d$.

Finally, we verify that $\mathcal{M}(D)$ satisfies the required operators and cokernel conditions of Definition \ref{dfn:m_functor}. For each $i \in \Delta$, $\mathcal{M}(D)$ is stable under $\varphi_i$ and $N_{\nabla}^i$ since each $\mathcal{M}_{i, \Delta}(D)$ is, and by Lemma \ref{lemma:m_i_free} there exists an integer $h_i \ge 0$ such that the cokernel of the linearized partial Frobenius $1 \otimes \varphi_i$ on $\mathcal{M}_i(D)$ (and hence on $\mathcal{M}(D)$ by localizing at the height-1 primes dividing $E_i(u_i)$) is annihilated by $E_i(u_i)^{h_i}$. Since the partial Frobenius maps commute, the linearization of the total Frobenius map $1 \otimes \varphi_\Delta$ factors as the composition of the linearizations $1 \otimes \varphi_i$, so its cokernel is annihilated by the product $E_\Delta^{\underline{h}}(\underline{u}) := \prod_{i \in \Delta} E_i(u_i)^{h_i}$:
\[
E_\Delta^{\underline{h}}(\underline{u}) \cdot \mathrm{coker}(1 \otimes \varphi_\Delta) = 0.
\]
\end{proof}

\subsection*{The $\mathcal{D}$ functor}
The space $\mathcal{O}_\Delta$ is a Fréchet space, whose topology is defined by the family of seminorms 
\[
    \|f\|_{\underline{r}} := \sup_{x\in D([0, r_1)) \times \dots \times D([0, r_\delta))} |f(x)|
\]
for each multi-index $\underline{r}:=(r_1, \dots, r_\delta) \in (0, 1)^\delta$. Furthermore, on any finite free module over $\mathcal{O}_\Delta$, we can define corresponding seminorms for each $\underline{r} \in (0, 1)^\delta$ by taking the maximum of the seminorms $\| \cdot \|_{\underline{r}}$ of the coordinates of an element. The norm of a subset of the module is then defined as the supremum of the norms of its elements. This topological framework will be essential in the subsequent lemma, where we construct a map via a telescoping sum.

\begin{lemma}\label{lemma:xi}
Let $M \in \mathrm{Mod}_{\mathcal{O}_\Delta}^{(\varphi_i, N_\nabla^i)_{i\in \Delta}}$.
\begin{enumerate}
    \item For each $i \in \Delta$, there is a unique $\mathcal{O}_\Delta$-linear, $(\varphi_j)_{j\in\Delta}$-equivariant morphism
    \[
    \xi_i : \mathcal{D}_i(M) \otimes_{\mathcal{O}_\Delta^{(i)}} \mathcal{O}_\Delta \longrightarrow M
    \]
    whose reduction modulo $u_i$ induces the identity on $\mathcal{D}_i(M) = M/u_i M$. The map $\xi_i$ is injective, its cokernel is annihilated by $\lambda_i^{h_i}$, and for any $r_i \in (|\pi_i|, |\pi_i|^{1/p})$, the image of $\xi_i$ on $D([0, r_i)) \times D([0,1)^{\delta-1})$ coincides with the image of $1 \otimes \varphi_{M,i} : \varphi_i^* M \to M$.
    \item Composing the directional morphisms $\xi_i$ yields the unique $\mathcal{O}_\Delta$-linear, $(\varphi_i)_{i\in\Delta}$-equivariant (and hence $\varphi_\Delta$-equivariant) morphism
    \[
    \xi : \mathcal{D}(M) \otimes_{K_{0,\Delta}} \mathcal{O}_\Delta \longrightarrow M
    \]
    whose reduction modulo $I := (u_1, \dots, u_\delta)$ induces the identity on $\mathcal{D}(M)$. Furthermore, $\xi$ is injective, its cokernel is annihilated by a finite power of $\lambda_\Delta$, and for $r_i \in (|\pi_i|, |\pi_i|^{1/p})$ ($i \in \Delta$), the image of $\xi$ on $D([0, r_1) \times \dots \times [0, r_\delta))$ coincides with the image of $1 \otimes \varphi_{M,\Delta} : \varphi_\Delta^* M \to M$.
\end{enumerate}
\end{lemma}

\begin{proof}
(1) Fix $i \in \Delta$, and set $\mathcal{O}_\Delta^{(i)} := \mathcal{O}_\Delta / u_i \mathcal{O}_\Delta \cong \mathcal{O}_{\Delta \setminus \{i\}} \otimes_{K_{0,\Delta \setminus \{i\}}} K_{0,\Delta}$. Since $M$ is a finite free $\mathcal{O}_\Delta$-module, $\mathcal{D}_i(M) = M/u_i M$ is a finite free $\mathcal{O}_\Delta^{(i)}$-module on which $\varphi_{\mathcal{D}_i(M),i}$ is $\varphi_i$-semilinear and bijective (as $E_i(0) \in K_{0,i}^\times$ is invertible). Choose an $\mathcal{O}_\Delta^{(i)}$-linear section $s_{0,i} \in \mathrm{Hom}_{\mathcal{O}_\Delta^{(i)}}(\mathcal{D}_i(M), M)$ reducing to the identity modulo $u_i$. For each integer $m_i \ge 0$, define the $\mathcal{O}_\Delta^{(i)}$-linear approximation
\[
s_{m_i}^{(i)} := \varphi_{M,i}^{m_i} \circ s_{0,i} \circ \varphi_{\mathcal{D}_i(M),i}^{-m_i} \in \mathrm{Hom}_{\mathcal{O}_\Delta^{(i)}}(\mathcal{D}_i(M), M).
\]
Because $s_{0,i} \equiv \mathrm{id}_{\mathcal{D}_i(M)} \pmod{u_i}$ and $\varphi_{M,i} \bmod u_i = \varphi_{\mathcal{D}_i(M),i}$, the single-step difference $t_i := \varphi_{M,i} \circ s_{0,i} \circ \varphi_{\mathcal{D}_i(M),i}^{-1} - s_{0,i}$ has image contained in $u_i M$. Since $\varphi_i^{m_i}(u_i) = u_i^{p^{m_i}}$, we have
\[
s_{m_i+1}^{(i)} - s_{m_i}^{(i)} = u_i^{p^{m_i}} \varphi_{M,i}^{m_i} \circ \left( u_i^{-1} t_i \right) \circ \varphi_{\mathcal{D}_i(M),i}^{-m_i}.
\]
On any polydisk of radius $\underline{r} = (r_1, \dots, r_\delta) \in (0,1)^\delta$, there exists $C_i > 0$ such that $\|\varphi_{M,i}\|_{\underline{r}} \le p^{C_i}$ and $\|\varphi_{\mathcal{D}_i(M),i}^{-1}\|_{\underline{r}} \le p^{C_i}$. Thus $\|s_{m_i+1}^{(i)} - s_{m_i}^{(i)}\|_{\underline{r}} \le r_i^{p^{m_i}} p^{2C_i m_i} \|u_i^{-1} t_i\|_{\underline{r}} \to 0$ as $m_i \to \infty$. By completeness of $M$ in the Fr\'echet topology, $s_{m_i}^{(i)}$ converges to a unique $\varphi_i$-equivariant $\mathcal{O}_\Delta^{(i)}$-linear map $\xi_{0,i} : \mathcal{D}_i(M) \to M$ lifting $\mathrm{id}_{\mathcal{D}_i(M)}$. Since $\varphi_j$ ($j \ne i$) commutes with $\varphi_i$ and acts as the identity on $u_i$, uniqueness of the $\varphi_i$-equivariant lift implies that $\xi_{0,i}$ also commutes with $\varphi_j$ for all $j \ne i$. Extending by $\mathcal{O}_\Delta$-linearity gives $\xi_i : \mathcal{D}_i(M) \otimes_{\mathcal{O}_\Delta^{(i)}} \mathcal{O}_\Delta \to M$.

To see that $\xi_i$ is injective with cokernel killed by $\lambda_i^{h_i}$, choose an $\mathcal{O}_\Delta$-linear map $s_{0,i}' : M \to \mathcal{D}_i(M) \otimes_{\mathcal{O}_\Delta^{(i)}} \mathcal{O}_\Delta$ reducing to the identity modulo $u_i$. Since $E_i(u_i)^{h_i} \mathrm{coker}(1 \otimes \varphi_{M,i}) = 0$, the reverse sequence
\[
    s_{m_{i}}^{\prime\,(i)} := (1 \otimes \varphi_{\mathcal{D}_{i}(M),i})^{m_{i}} \circ \varphi_i^{* m_{i}}(s_{0,i}') \circ (1 \otimes \varphi_{M,i})^{-m_{i}}
\]
has denominators dividing $\prod_{k=0}^{m_i-1} \varphi_i^k(E_i(u_i))^{h_i} \mid \lambda_i^{h_i}$, and $\lambda_i^{h_i} s_{m_i}^{\prime\,(i)}$ converges by the same estimate to an $\mathcal{O}_\Delta$-linear map $\xi_i' : M \to \mathcal{D}_i(M) \otimes_{\mathcal{O}_{\Delta^{(i)}}} \lambda_i^{-h_i}\mathcal{O}_\Delta$ satisfying $\xi_i' \circ \xi_i = 1$ and $\xi_i \circ \xi_i' = 1$ after inverting $\lambda_i$. If $r_i \in (|\pi_i|, |\pi_i|^{1/p})$, then $\varphi_i(\lambda_i)$ is a unit on $D([0, r_i)) \times D([0,1)^{\delta-1})$, so $\varphi_i^*(\xi_i)$ is an isomorphism there, and $\xi_i \circ (1 \otimes \varphi_{\mathcal{D}_i(M),i}) = (1 \otimes \varphi_{M,i}) \circ \varphi_i^*(\xi_i)$ implies that the image of $\xi_i$ coincides with the image of $1 \otimes \varphi_{M,i}$.

(2) Ordering $\Delta = \{1, \dots, \delta\}$ and setting $M^{(i)} := M/(u_1, \dots, u_i)M = \mathcal{D}_i(M^{(i-1)})$ (with $M^{(0)} = M$ and $M^{(\delta)} = \mathcal{D}(M)$), part (1) applied inductively to $M^{(i-1)}$ in the variable $u_i$ produces $\mathcal{O}_\Delta$-linear, $(\varphi_j)_{j\in\Delta}$-equivariant morphisms
\[
\xi^{(i)} : M^{(i)} \otimes_{\mathcal{O}_{\Delta \setminus \{1,\dots,i\}}} \mathcal{O}_\Delta \longrightarrow M^{(i-1)} \otimes_{\mathcal{O}_{\Delta \setminus \{1,\dots,i-1\}}} \mathcal{O}_\Delta.
\]
Defining $\xi := \xi^{(1)} \circ \xi^{(2)} \circ \dots \circ \xi^{(\delta)} : \mathcal{D}(M) \otimes_{K_{0,\Delta}} \mathcal{O}_\Delta \to M$ (which is the limit of the composite approximations $s_{m_1,\dots,m_\delta} = s_{m_1}^{(1)} \circ \dots \circ s_{m_\delta}^{(\delta)} = \left(\prod_{i\in\Delta} \varphi_{M,i}^{m_i}\right) \circ s_0 \circ \left(\prod_{i\in\Delta} \varphi_{\mathcal{D}(M),i}^{-m_i}\right)$ as $m_1, \dots, m_\delta \to \infty$ sequentially), we obtain a $(\varphi_i)_{i\in\Delta}$-equivariant morphism reducing to $\mathrm{id}_{\mathcal{D}(M)}$ modulo $I$. Since each $\xi^{(i)}$ is injective with cokernel killed by $\lambda_i^{h_i}$, $\xi$ is injective with cokernel killed by $\prod_{i\in\Delta} \lambda_i^{h_i} \mid \lambda_\Delta^h$ ($h = \max_i h_i$), and its image on $D([0, r_1) \times \dots \times [0, r_\delta))$ coincides with the image of $1 \otimes \varphi_{M,\Delta}$.
\end{proof}

\begin{pro}\label{pro:filtration}
Let $M \in Mod_{\mathcal{O}_\Delta}^{(\varphi_i, N_\nabla^{i})_{i\in \Delta}}$, then $$\mathcal{D}(M) \in  MF_{K_\Delta}^{(\varphi_i, N_i)_{i\in \Delta}, \text{ Fil}\geq 0}.$$
\end{pro}

\begin{proof}
For every $i \in \Delta$, we define $\mathcal{D}_i(M)$ as the module $M$ modulo $u_i$, so $\mathcal{D}_i(M)$ is a $\mathcal{O}_\Delta^{(i)}(:=K_i\otimes_{K_{0, \Delta}} \mathcal{O}_\Delta)$-module, on which a partial Frobenius map $\varphi_i$ acting on $K_i$ and a monodromy operator $N_i$ on $K_i$ are uniquely defined by modulo $u_i$. Furthermore, for the other $j\neq i$ partial Frobenii $\varphi_j$ on $\mathcal{O}_\Delta^{(i)}$ and the partial differential operators $N_\nabla^{j}$ on $\mathcal{O}_\Delta^{(i)}$ remain. The commutativity remains in the case $j \neq i$ and for $i$ it follows that $p \cdot ( \varphi_i \circ N_i) = N_i \circ \varphi_i$. 

The final object $\mathcal{D}(M)$ is independent of the order in which these univariate functors are applied. By standard module isomorphism theorems, we have, for example,
\[
    (M / u_1 M) / u_2 (M / u_1 M) \cong M / (u_1, u_2) M \cong (M / u_2 M) / u_1 (M / u_2 M).
\]
Because the operators commute in the original space, their induced actions on the successive quotients necessarily commute as well. Consequently, $\mathcal{D}(M)$ is a finite free $K_\Delta$-module equipped with the appropriate partial Frobenii and monodromy operators.

The multifiltration on $\mathcal{D}(M)_{K_\Delta} := \mathcal{D}(M) \otimes_{K_{0,\Delta}} K_\Delta$ is constructed from the directional filtrations on $\mathcal{D}(M)_{K_i} := \mathcal{D}(M) \otimes_{K_{0,i}} K_i$ ($i \in \Delta$) via the transversality of the functors $\mathcal{D}_i$ and the commutativity of the partial operators. Specifically, for each $i \in \Delta$, let $p_i : \mathcal{D}_i(M) = M/u_i M \twoheadrightarrow \mathcal{D}(M) = M/(u_1, \dots, u_\delta)M$ denote the natural projection modulo $(u_j)_{j \ne i}$.

To construct the filtration on $\mathcal{D}_i(M)_{K_i} := \mathcal{D}_i(M) \otimes_{K_{0,i}} K_i$ (and hence on $\mathcal{D}(M)_{K_i}$), we first define a decreasing, separated, and exhaustive filtration on $\varphi_i^* M$ by powers of the Eisenstein polynomial $E_i(u_i)$:
\[
\mathrm{Fil}^{k_i}(\varphi_i^* M) := \left\{ x \in \varphi_i^* M \;\middle|\; (1 \otimes \varphi_{M,i})(x) \in E_i(u_i)^{k_i} M \right\}, \quad k_i \in \mathbb{Z}.
\]
Equipping $(1 \otimes \varphi_{M,i})(\varphi_i^* M) \subset M$ with the image filtration $\mathrm{Fil}^{k_i}((1 \otimes \varphi_{M,i})(\varphi_i^* M)) := (1 \otimes \varphi_{M,i})(\varphi_i^* M) \cap E_i(u_i)^{k_i} M$, we restrict to the cylinder $D([0, r_i)) \times D([0, 1)^{\delta-1})$ for any fixed $r_i \in (|\pi_i|, |\pi_i|^{1/p})$. By Lemma \ref{lemma:xi}(1), the morphism $\xi_i : \mathcal{D}_i(M) \otimes_{\mathcal{O}_\Delta^{(i)}} \mathcal{O}_\Delta \to M$ maps $(\mathcal{D}_i(M) \otimes_{\mathcal{O}_\Delta^{(i)}} \mathcal{O}_\Delta)_{[0, r_i)}$ isomorphically onto $(1 \otimes \varphi_{M,i})(\varphi_i^* M)_{[0, r_i)}$. Pulling back the filtration along $\xi_i$ and reducing modulo $E_i(u_i)$ via the natural quotient map
\[
\left(\mathcal{D}_i(M) \otimes_{\mathcal{O}_\Delta^{(i)}} \mathcal{O}_\Delta\right)_{[0, r_i)} \twoheadrightarrow \left(\mathcal{D}_i(M) \otimes_{\mathcal{O}_\Delta^{(i)}} \mathcal{O}_\Delta\right)_{[0, r_i)} \big/ E_i(u_i) \xrightarrow{\;\sim\;} \mathcal{D}_i(M) \otimes_{K_{0,i}} K_i = \mathcal{D}_i(M)_{K_i},
\]
we obtain a decreasing, separated, and exhaustive filtration $\mathrm{Fil}^{k_i}(\mathcal{D}_i(M)_{K_i})$ by $\mathcal{O}_\Delta^{(i)} \otimes_{K_{0,i}} K_i$-submodules, which is effective ($\mathrm{Fil}^0 = \mathcal{D}_i(M)_{K_i}$) since $(1 \otimes \varphi_{M,i})(\varphi_i^* M) \subset M$. Applying the projection $p_i \otimes \mathrm{id}_{K_i} : \mathcal{D}_i(M)_{K_i} \twoheadrightarrow \mathcal{D}(M)_{K_i}$ (i.e., reducing modulo $u_j$ for all $j \ne i$) defines a decreasing, separated, exhaustive, and effective filtration $\mathrm{Fil}^{k_i}(\mathcal{D}(M)_{K_i}) := (p_i \otimes \mathrm{id}_{K_i})(\mathrm{Fil}^{k_i}(\mathcal{D}_i(M)_{K_i}))$ on $\mathcal{D}(M)_{K_i}$. Finally, for each $\underline{k} = (k_1, \dots, k_\delta) \in \mathbb{Z}^\delta$, we define the multifiltration on $\mathcal{D}(M)_{K_\Delta}$ as in Definition \ref{dfn:phi-Nmod}(3) by
\[
\mathrm{Fil}^{\underline{k}}(\mathcal{D}(M)_{K_\Delta}) := \bigcap_{i \in \Delta} \left( \mathrm{Fil}^{k_i}(\mathcal{D}(M)_{K_i}) \otimes_{(K_i \otimes_{K_{0,i}} K_{0,\Delta})} K_\Delta \right) \subset \mathcal{D}(M)_{K_\Delta}.
\]
\end{proof}

\subsection*{The natural isomorphisms}
To establish the equivalence of categories, we construct the natural isomorphisms $\mathcal{D}(\mathcal{M}(D)) \simeq D$ and $\mathcal{M}(\mathcal{D}(M)) \simeq M$. The first direction is relatively straightforward: the transversality of the partial modules $\mathcal{O}_\Delta \otimes_{\mathcal{O}_i \otimes K_{0, \Delta}} \mathcal{M}_i(D)$ allows us to safely interchange the order of the intersection and the univariate reductions $\mathcal{D}_i$.

\begin{pro}\label{pro:first}
Let $D \in \mathrm{MF}_{K_{0,\Delta}}^{(\varphi_i, N_i)_{i \in \Delta}, \mathrm{Fil} \ge 0}$, then there exists a natural isomorphism
\[
\mathcal{D}(\mathcal{M}(D)) \xrightarrow{\;\sim\;} D
\]
in the category $\mathrm{MF}_{K_{0,\Delta}}^{(\varphi_i, N_i)_{i \in \Delta}, \mathrm{Fil} \ge 0}$.
\end{pro}

\begin{proof}
Recall that $\mathcal{D}(M) = \mathcal{D}_\delta \circ \dots \circ \mathcal{D}_1(M) \cong M / I M$, where $I := (u_1, \dots, u_\delta) \subset \mathcal{O}_\Delta$, and that each $\mathcal{M}_{i,\Delta}(D) = \mathcal{O}_\Delta[1/\lambda_{\Delta \setminus \{i\}}] \otimes_{(\mathcal{O}_i \otimes_{K_{0,i}} K_{0,\Delta})} \mathcal{M}_i(D)$ is a submodule of the ambient space $\mathcal{D}_0[1/\lambda_\Delta]$, where
\[
\mathcal{D}_0 = (\mathcal{O}_\Delta[l_{u_1}, \dots, l_{u_\delta}] \otimes_{K_{0,\Delta}} D)^{N_1 = \dots = N_\delta = 0}.
\]
For each $i \in \Delta$, the infinite product $\lambda_i = \prod_{n=0}^\infty \varphi_i^n(E_i(u_i)/E_i(0)) \in \mathcal{O}_i$ satisfies $\lambda_i(0) = 1$, and is therefore an invertible unit in the local ring $\mathcal{O}_{\Delta, I}$ at the origin $u_1 = \dots = u_\delta = 0$. Consequently, at the origin $I = (u_1, \dots, u_\delta)$, all denominators $\lambda_1, \dots, \lambda_\delta$ are units, so the inclusions
\[
\mathcal{D}_0 \subseteq \mathcal{M}(D) = \bigcap_{i \in \Delta} \mathcal{M}_{i,\Delta}(D) \subseteq \mathcal{D}_0[1/\lambda_\Delta]
\]
become equalities after localizing at $I$. In particular, reducing the finite intersection $\mathcal{M}(D) = \bigcap_{i \in \Delta} \mathcal{M}_{i,\Delta}(D)$ modulo $I = (u_1, \dots, u_\delta)$ commutes with the intersection and yields:
\begin{align*}
\mathcal{D}(\mathcal{M}(D)) &= \mathcal{D}\!\left(\bigcap_{i \in \Delta} \mathcal{M}_{i,\Delta}(D)\right) = \bigcap_{i \in \Delta} \mathcal{D}(\mathcal{M}_{i,\Delta}(D)) = \mathcal{D}(\mathcal{D}_0) \\
&= (K_{0,\Delta}[l_{u_1}, \dots, l_{u_\delta}] \otimes_{K_{0,\Delta}} D)^{N_1 = \dots = N_\delta = 0} \subset K_{0,\Delta}[l_{u_1}, \dots, l_{u_\delta}] \otimes_{K_{0,\Delta}} D.
\end{align*}
In Lemma \ref{lemma:eta}, we showed that the specialization map $\eta : (l_{u_1}, \dots, l_{u_\delta}) \mapsto (0, \dots, 0)$ defines a natural isomorphism of $(\varphi_i, N_i)_{i \in \Delta}$-modules over $K_{0,\Delta}$:
\[
\eta : \mathcal{D}(\mathcal{M}(D)) = (K_{0,\Delta}[l_{u_1}, \dots, l_{u_\delta}] \otimes_{K_{0,\Delta}} D)^{N_1 = \dots = N_\delta = 0} \xrightarrow{\;\sim\;} D.
\]
Furthermore, since for each $i \in \Delta$ the univariate reduction $\mathcal{D}_i(\mathcal{M}_i(D)) \xrightarrow{\sim} D$ strictly preserves $\varphi_i$, $N_i$, and the directional filtration $\mathrm{Fil}^{k_i}(D_{K_i})$ by \cite[Proposition 1.2.8]{Kis1}, and the multifiltration on $\mathcal{D}(\mathcal{M}(D))_{K_\Delta}$ is defined as the intersection of the base-changed directional filtrations (Proposition \ref{pro:filtration}), $\eta$ is a strict isomorphism of filtered $(\varphi_i, N_i)_{i \in \Delta}$-modules.
\end{proof}

For a module $M \in \text{Mod}_{\mathcal{O}_\Delta}^{(\varphi_i, N_\nabla^{i})_{i\in \Delta}}$, we define 
\[
    \mathcal{D}_0(M) := (\mathcal{O}_\Delta[l_{u_1}, \dots, l_{u_\delta}] \otimes_{K_{0,\Delta}} \mathcal{D}(M))^{N_1=\dots=N_\delta=0}.
\]
The filtration on $\mathcal{D}_0(M)$ is inherited from the filtration on $\mathcal{D}(M)$, as defined in Proposition \ref{pro:filtration}. Using the maps $\eta$ and $\xi$ defined in Lemma \ref{lemma:eta} and Lemma \ref{lemma:xi}, respectively, we obtain the following composition of maps from $\mathcal{D}_0(M)$ to $M$:
\[
    \mathcal{D}_0(M) = (K_{0,\Delta}[l_{u_1},\dots, l_{u_\delta}]\otimes_{K_{0, \Delta}} \mathcal{D}(M))^{N_1=\dots=N_\delta=0}\otimes_{K_{0, \Delta}} \mathcal{O}_\Delta \xrightarrow{\eta\otimes 1} \mathcal{D}(M)\otimes_{K_{0, \Delta}} \mathcal{O}_\Delta \xrightarrow{\xi} M.
\]

\begin{pro}\label{pro:second}
Let $M \in \mathrm{Mod}_{\mathcal{O}_\Delta}^{(\varphi_i, N_\nabla^{i})_{i\in \Delta}}$, then there exists a natural isomorphism 
\[
    \mathcal{M}(\mathcal{D}(M))\xrightarrow{\;\sim\;} M
\]
in the category $\mathrm{Mod}_{\mathcal{O}_\Delta}^{(\varphi_i, N_\nabla^{i})_{i\in \Delta}}$.
\end{pro}

\begin{proof}
Let $M':=\mathcal{M}(\mathcal{D}(M))$ and let $\mathcal{D}_0(M)$ be as before. Via the injective map $\xi \circ (\eta \otimes 1) : \mathcal{D}_0(M) \to M$ whose cokernel is killed by a finite power of $\lambda_\Delta$, we identify both $M$ and $M'$ with $\mathcal{O}_\Delta$-submodules of $\mathcal{D}_0(M)[1/\lambda_\Delta]$. Since $M$ and $M'$ both have finite $E_\Delta$-height, it suffices to prove that $M = M'$ over the polydisk $D([0, r)^\delta)$, where $r=(r_1,\ldots,r_\delta)$ with $r_i \in (|\pi_i|, |\pi_i|^{1/p})$ for every $i\in \Delta$, as this equality then extends to the entire polydisk $D([0,1)^\delta)$ by pulling back along powers of the total Frobenius $\varphi_\Delta$. Away from the union of the divisors $E_i(u_i) = 0$ on $D([0, r)^\delta)$, all $\lambda_i$ are invertible, so $M = \mathcal{D}_0(M) = M'$ outside the hyperplanes $E_i(u_i) = 0$.

Write $D([0, r)^\delta) = \bigcup_{k\geq 0} U_{\rho_k}$, where $U_{\rho_k}:=\{x \in D([0, r)^\delta) \mid \forall j\in \Delta, \; |x_j|\leq \rho_{k,j} < r_j \}$ with $U_{\rho_{k}} \subseteq U_{\rho_{k+1}}$ and $\lim_{k\rightarrow \infty} \rho_{k, j} = r_j$ for every $j\in \Delta$. To prove that $M$ and $M'$ coincide as global $\mathcal{O}_\Delta$-modules, we rely on their structure as coherent sheaves on the quasi-Stein space $D([0, r)^\delta)$, prove that $M(U_{\rho_k}) = M'(U_{\rho_k})$ on each affinoid subdomain $U_{\rho_k}$, and take global sections.

By Lemma \ref{lemma:components2}, the coordinate ring $\mathcal{O}_\Delta(U_{\rho_k}) \cong \prod_{\alpha \in \Sigma} e_\alpha \mathcal{O}_\Delta(U_{\rho_k})$ of the closed affinoid polydisk $U_{\rho_k}$ is a finite direct product of free Tate algebras over the fields $F_\alpha = e_\alpha K_{0,\Delta}$, each of which is a regular Noetherian normal domain. Since $M(U_{\rho_k})$ and $M'(U_{\rho_k})$ are free (hence reflexive) $\mathcal{O}_\Delta(U_{\rho_k})$-submodules of $\mathcal{D}_0(M)(U_{\rho_k})[1/\lambda_\Delta]$, the discussion following Lemma \ref{lemma:reflexiv} implies that
\[
M(U_{\rho_k}) = \bigcap_{\mathrm{ht}(\mathfrak{p})=1} M(U_{\rho_k})_{\mathfrak{p}} \quad \text{and} \quad M'(U_{\rho_k}) = \bigcap_{\mathrm{ht}(\mathfrak{p})=1} M'(U_{\rho_k})_{\mathfrak{p}},
\]
where the intersections run over all height-$1$ prime ideals $\mathfrak{p}$ of $\mathcal{O}_\Delta(U_{\rho_k})$. Thus, to prove $M(U_{\rho_k}) = M'(U_{\rho_k})$, it suffices to show that $M(U_{\rho_k})_{\mathfrak{p}} = M'(U_{\rho_k})_{\mathfrak{p}}$ at every height-$1$ prime ideal $\mathfrak{p}$ of $\mathcal{O}_\Delta(U_{\rho_k})$.

Notice that on $U_{\rho_k}$ (where $\rho_{k,j} < r_j < |\pi_j|^{1/p}$), the factors $\varphi_j^n(E_j(u_j))$ for $n \ge 1$ have no zeros and are therefore units in $\mathcal{O}_\Delta(U_{\rho_k})$. Thus, for any height-$1$ prime ideal $\mathfrak{p}$ of $\mathcal{O}_\Delta(U_{\rho_k})$ that does not contain $E_i(u_i)$ for any $i \in \Delta$, all $\lambda_i$ are units in $\mathcal{O}_\Delta(U_{\rho_k})_{\mathfrak{p}}$, so $M(U_{\rho_k})_{\mathfrak{p}} = \mathcal{D}_0(M)(U_{\rho_k})_{\mathfrak{p}} = M'(U_{\rho_k})_{\mathfrak{p}}$.

It remains to check the height-$1$ prime ideals $\mathfrak{p}$ containing $E_i(u_i)$ for some $i \in \Delta$. By symmetry, it suffices to consider $i = 1$. Since $E_1(u_1) \in K_{0,1}[u_1]$ is Eisenstein over $K_{0,1}$ and each field factor $F_\alpha / K_{0,1}$ of $K_{0,\Delta}$ is unramified, $E_1(u_1)$ remains irreducible in each $F_\alpha[u_1]$ and generates a height-$1$ prime ideal in each connected component $e_\alpha \mathcal{O}_\Delta(U_{\rho_k})$. Let $\mathfrak{p}$ be any height-$1$ prime ideal of $\mathcal{O}_\Delta(U_{\rho_k})$ containing $E_1(u_1)$ (corresponding to a unique connected component $\alpha \in \Sigma$), and let $\mathcal{O}_{\mathfrak{p}} := \mathcal{O}_\Delta(U_{\rho_k})_{\mathfrak{p}} \cong (e_\alpha \mathcal{O}_\Delta(U_{\rho_k}))_{(E_1(u_1))}$ denote the localized ring, which is a discrete valuation ring with uniformizer $E_1(u_1)$. Let $M_{\mathfrak{p}} := M(U_{\rho_k})_{\mathfrak{p}}$ and $M'_{\mathfrak{p}} := M'(U_{\rho_k})_{\mathfrak{p}}$ be the localizations of $M(U_{\rho_k})$ and $M'(U_{\rho_k})$ at $\mathfrak{p}$. Both $M_{\mathfrak{p}}$ and $M'_{\mathfrak{p}}$ are finite free $\mathcal{O}_{\mathfrak{p}}$-modules inside $\mathcal{D}_0(M)(U_{\rho_k})_{\mathfrak{p}}[1/E_1(u_1)]$, and the cokernel of $1 \otimes \varphi_1$ on each is annihilated by $E_1(u_1)^{h_1}$ for a sufficiently large integer $h_1 \ge 0$.

Since $r_1 \in (|\pi_1|, |\pi_1|^{1/p})$, the elements $\varphi_1^n(E_1(u_1))$ for $n \ge 1$ are invertible on $U_{\rho_k}$, which implies that $\varphi_1(\lambda_1)$ is a unit in $\mathcal{O}_{\mathfrak{p}}$, because $M(U_{\rho_k})[1/\lambda_{\Delta}] = M'(U_{\rho_k})[1/\lambda_{\Delta}]$ and the other $\lambda_j$ ($j \neq 1$) are already units in $\mathcal{O}_{\mathfrak{p}}$, pulling back along $\varphi_1$ inverts $\varphi_1(\lambda_1)$ in $\mathcal{O}_{\mathfrak{p}}$ and yields an equality of the linearized Frobenius images:
\[
N_{\mathfrak{p}} := (1 \otimes \varphi_1)(\varphi_1^* M_{\mathfrak{p}}) = (1 \otimes \varphi_1)(\varphi_1^* M'_{\mathfrak{p}}) \subset M_{\mathfrak{p}} \cap M'_{\mathfrak{p}}.
\]
Equivalently, this follows from Lemma \ref{lemma:xi} and Proposition \ref{pro:first}, as both $(1 \otimes \varphi_1)(\varphi_1^* M_{\mathfrak{p}})$ and $(1 \otimes \varphi_1)(\varphi_1^* M'_{\mathfrak{p}})$ coincide with the image of $\mathcal{D}_1(M)(U_{\rho_k})_{\mathfrak{p}}$ under $\xi_1$.

We define two decreasing filtrations on the common free $\mathcal{O}_{\mathfrak{p}}$-module $N_{\mathfrak{p}}$ by setting, for each $m \in \mathbb{Z}$,
\[
\mathrm{Fil}^m_M N_{\mathfrak{p}} := N_{\mathfrak{p}} \cap E_1(u_1)^m M_{\mathfrak{p}} \quad \text{and} \quad \mathrm{Fil}^m_{M'} N_{\mathfrak{p}} := N_{\mathfrak{p}} \cap E_1(u_1)^m M'_{\mathfrak{p}}.
\]

\vspace{0.5cm}

\begin{lemma}
For every $m \in \mathbb{Z}$, the filtrations $\mathrm{Fil}^m_M N_{\mathfrak{p}}$ and $\mathrm{Fil}^m_{M'} N_{\mathfrak{p}}$ on $N_{\mathfrak{p}}$ coincide. Consequently, $M_{\mathfrak{p}} = M'_{\mathfrak{p}}$.
\end{lemma}

\vspace{0.5cm}

\begin{proof}
Let $\widehat{\mathcal{O}}_{\mathfrak{p}}$ denote the $E_1(u_1)$-adic completion of the discrete valuation ring $\mathcal{O}_{\mathfrak{p}}$, and let $\widehat{N}_{\mathfrak{p}} := \widehat{\mathcal{O}}_{\mathfrak{p}} \otimes_{\mathcal{O}_{\mathfrak{p}}} N_{\mathfrak{p}}$. Note that $\widehat{\mathcal{O}}_{\mathfrak{p}}$ naturally contains $e_\alpha \widehat{\mathfrak{S}}_{0,1}$ as a subring with the same uniformizer $E_1(u_1)$ (or $u_1 - \pi_1$), and $\widehat{\mathcal{O}}_{\mathfrak{p}}$ is faithfully flat over $\mathcal{O}_{\mathfrak{p}}$. By Lemma \ref{lemma:xi}, the univariate section $\xi_1$ identifies $N_{\mathfrak{p}}$ with $(\mathcal{D}_1(M) \otimes_{\mathcal{O}_\Delta^{(1)}} \mathcal{O}_\Delta)(U_{\rho_k})_{\mathfrak{p}}$, so passing to the completion identifies $\widehat{N}_{\mathfrak{p}}[1/E_1(u_1)]$ with $\widehat{\mathcal{O}}_{\mathfrak{p}}[1/E_1(u_1)] \otimes_{\mathcal{O}_\Delta^{(1)}} \mathcal{D}_1(M)$.

By the construction of the partial filtration in Proposition \ref{pro:filtration}, the filtration on $\mathcal{D}_1(M)_{K_1}$ is defined as the image of $\mathrm{Fil}^m_M N_{\mathfrak{p}}$ modulo $E_1(u_1)$. More precisely, if we equip $\widehat{\mathcal{O}}_{\mathfrak{p}}[1/E_1(u_1)]$ with its $E_1(u_1)$-adic filtration $\mathrm{Fil}^i \widehat{\mathcal{O}}_{\mathfrak{p}}[1/E_1(u_1)] = E_1(u_1)^i \widehat{\mathcal{O}}_{\mathfrak{p}}$ and $\mathcal{D}_1(M)_{K_1}$ with the filtration from Proposition \ref{pro:filtration} (and Kisin's comparison isomorphism \cite[Lemma 1.2.12(4)]{Kis1} along the $u_1$-direction), then the completed filtration $\widehat{\mathrm{Fil}}^m_M \widehat{N}_{\mathfrak{p}} := \widehat{\mathcal{O}}_{\mathfrak{p}} \otimes_{\mathcal{O}_{\mathfrak{p}}} \mathrm{Fil}^m_M N_{\mathfrak{p}}$ is recovered as
\[
\widehat{\mathrm{Fil}}^m_M \widehat{N}_{\mathfrak{p}} = \widehat{N}_{\mathfrak{p}} \cap \mathrm{Fil}^m\left(\widehat{\mathcal{O}}_{\mathfrak{p}}[1/E_1(u_1)] \otimes_{K_1} \mathcal{D}_1(M)_{K_1}\right).
\]
On the other hand, by Proposition \ref{pro:first}, the natural isomorphism $\mathcal{D}(M) \cong \mathcal{D}(\mathcal{M}(\mathcal{D}(M))) = \mathcal{D}(M')$ is a strict isomorphism of filtered $(\varphi_i, N_i)_{i \in \Delta}$-modules, meaning that $\mathcal{D}_1(M)_{K_1}$ and $\mathcal{D}_1(M')_{K_1}$ carry the exact same filtration. Since both $\widehat{\mathrm{Fil}}^m_M \widehat{N}_{\mathfrak{p}}$ and $\widehat{\mathrm{Fil}}^m_{M'} \widehat{N}_{\mathfrak{p}}$ are given by intersecting $\widehat{N}_{\mathfrak{p}}$ with the $m$-th step of the same filtration on $\widehat{\mathcal{O}}_{\mathfrak{p}}[1/E_1(u_1)] \otimes_{K_1} \mathcal{D}_1(M)_{K_1}$, we obtain
\[
\widehat{\mathrm{Fil}}^m_M \widehat{N}_{\mathfrak{p}} = \widehat{\mathrm{Fil}}^m_{M'} \widehat{N}_{\mathfrak{p}} \quad \text{for all } m \in \mathbb{Z}.
\]
Since $\mathcal{O}_{\mathfrak{p}} \to \widehat{\mathcal{O}}_{\mathfrak{p}}$ is faithfully flat, intersecting with $N_{\mathfrak{p}}$ inside $\widehat{N}_{\mathfrak{p}}$ yields $\mathrm{Fil}^m_M N_{\mathfrak{p}} = \mathrm{Fil}^m_{M'} N_{\mathfrak{p}}$ for all $m \in \mathbb{Z}$.

Finally, since the cokernel of $1 \otimes \varphi_1$ on both $M_{\mathfrak{p}}$ and $M'_{\mathfrak{p}}$ is annihilated by $E_1(u_1)^{h_1}$, we have $E_1(u_1)^{h_1} M_{\mathfrak{p}} \subset N_{\mathfrak{p}}$ and $E_1(u_1)^{h_1} M'_{\mathfrak{p}} \subset N_{\mathfrak{p}}$. By definition of the filtrations at $m = h_1$, this implies
\[
E_1(u_1)^{h_1} M_{\mathfrak{p}} = N_{\mathfrak{p}} \cap E_1(u_1)^{h_1} M_{\mathfrak{p}} = \mathrm{Fil}^{h_1}_M N_{\mathfrak{p}} = \mathrm{Fil}^{h_1}_{M'} N_{\mathfrak{p}} = N_{\mathfrak{p}} \cap E_1(u_1)^{h_1} M'_{\mathfrak{p}} = E_1(u_1)^{h_1} M'_{\mathfrak{p}}.
\]
Multiplying by $E_1(u_1)^{-h_1}$ inside $N_{\mathfrak{p}}[1/E_1(u_1)]$ yields $M_{\mathfrak{p}} = M'_{\mathfrak{p}}$.
\end{proof}

\vspace{0.5cm}

By the lemma above, the localizations of $M(U_{\rho_k})$ and $M'(U_{\rho_k})$ coincide at every height-$1$ prime ideal $\mathfrak{p}$ dividing $E_1(u_1)$, and by symmetry at every height-$1$ prime ideal dividing $E_i(u_i)$ for all $i \in \Delta$. Since they also coincide at all other height-$1$ prime ideals of $\mathcal{O}_\Delta(U_{\rho_k})$, Lemma \ref{lemma:reflexiv} (applied componentwise via Lemma \ref{lemma:components2}) implies that $M(U_{\rho_k}) = M'(U_{\rho_k})$ inside $\mathcal{D}_0(M)(U_{\rho_k})[1/\lambda_\Delta]$.

Because the closed affinoid subdomains $U_{\rho_k}$ are nested and the identity inclusions on $U_{\rho_k}$ are compatible with restriction, they glue to an equality of coherent sheaves $M|_{D([0,r)^{\delta})} = M'|_{D([0,r)^{\delta})}$ on the quasi-Stein space $D([0,r)^{\delta})$. Extending via the total Frobenius $\varphi_{\Delta}$ to the entire polydisk $D([0,1)^{\delta})$ and taking global sections via Kiehl's Theorem B (Theorem \ref{thm:kiehl}) yields $M = M'$ as $\mathcal{O}_{\Delta}$-modules.
\end{proof}

\subsection*{End of the proof}
\begin{pro}\label{pro:exact_tensor}
$\mathcal{M}$ and $\mathcal{D}$ are exact, tensor functors between the categories $\text{MF}_{K_{0,\Delta}}^{(\varphi_i, N_i)_{i\in \Delta}, \text{ Fil}\geq 0}$ and $\text{Mod}_{\mathcal{O}_\Delta}^{(\varphi_i, N_\nabla^{i})_{i\in \Delta}}$.
\end{pro}
\begin{proof}
The exactness and tensor compatibility of the functor $\mathcal{D}$ follow immediately from its definition and the strict compatibility of the induced filtrations (as in \cite[Theorem 1.2.15]{Kis1}). The main difficulties lie with the functor $\mathcal{M}$. In the univariate case, the points where exactness and the tensor property are not a priori guaranteed correspond to the roots of $\lambda$, which form a discrete subset of the disk $D([0, 1))$. In the multivariable case, however, the zero locus of $\lambda$ forms a union of hyperplanes in the poly-disk $D([0,1))^\delta$.

To prove that $\mathcal{M}$ is exact, let
\[
0 \longrightarrow D' \longrightarrow D \xrightarrow{\;q\;} D'' \longrightarrow 0
\]
be a short exact sequence in $\mathrm{MF}_{K_{0,\Delta}}^{(\varphi_i, N_i)_{i \in \Delta}, \mathrm{Fil} \ge 0}$. We must show that the induced sequence of $\mathcal{O}_\Delta$-modules
\[
0 \longrightarrow \mathcal{M}(D') \longrightarrow \mathcal{M}(D) \xrightarrow{\;\mathcal{M}(q)\;} \mathcal{M}(D'') \longrightarrow 0
\]
is exact. Consider the exhausting covering of $D([0,1)^\delta) = \bigcup_{k \ge 1} U_{\rho_k}$ by nested closed affinoid polydisks $U_{\rho_k}$ with coordinate rings $A_k := \mathcal{O}_\Delta(U_{\rho_k}) = \prod_{\alpha \in \Sigma} e_\alpha A_k$ (Lemma \ref{lemma:components2}). For a fixed $k \ge 1$, set $\mathfrak{M}_{D'} := A_k \otimes_{\mathcal{O}_\Delta} \mathcal{M}(D')$, $\mathfrak{M}_D := A_k \otimes_{\mathcal{O}_\Delta} \mathcal{M}(D)$, and $\mathfrak{M}_{D''} := A_k \otimes_{\mathcal{O}_\Delta} \mathcal{M}(D'')$, which are finite free $A_k$-modules by Proposition \ref{pro:m_d_category}.

By construction of the functor $\mathcal{M}$, we have $\mathfrak{M}_{D'} = \mathfrak{M}_D \cap \mathcal{D}_{0,k}(D')[1/\lambda_\Delta]$ inside $\mathcal{D}_{0,k}(D)[1/\lambda_\Delta]$, which immediately implies that $\mathfrak{M}_{D'} \to \mathfrak{M}_D$ is injective and equals the kernel of the induced map $q_k : \mathfrak{M}_D \to \mathfrak{M}_{D''}$. To prove that $q_k : \mathfrak{M}_D \to \mathfrak{M}_{D''}$ is surjective over $A_k$, we first check surjectivity after localizing at any height-$1$ prime ideal $\mathfrak{p} \subset A_k$. By Lemma \ref{lem:height1_localization}, $A_{k, \mathfrak{p}}$ is a discrete valuation ring (lying on a unique connected component $e_\alpha A_k$), and at most one $\lambda_i$ belongs to $\mathfrak{p}$:

If $\lambda_j \notin \mathfrak{p}$ for all $j \in \Delta$, then by Lemma \ref{lem:height1_localization}(1) the sequence $0 \to (\mathfrak{M}_{D'})_{\mathfrak{p}} \to (\mathfrak{M}_D)_{\mathfrak{p}} \to (\mathfrak{M}_{D''})_{\mathfrak{p}} \to 0$ is identified with the base-changed sequence $0 \to \mathcal{D}_{0,k}(D')_{\mathfrak{p}} \to \mathcal{D}_{0,k}(D)_{\mathfrak{p}} \to \mathcal{D}_{0,k}(D'')_{\mathfrak{p}} \to 0$, which is exact by Lemma \ref{lemma:eta}.

If $\lambda_i \in \mathfrak{p}$ for a unique index $i \in \Delta$, then by Lemma \ref{lem:height1_localization}(2) we have $(\mathfrak{M}_{D'})_{\mathfrak{p}} = (A_k \otimes \mathcal{M}_i(D'))_{\mathfrak{p}}$, $(\mathfrak{M}_D)_{\mathfrak{p}} = (A_k \otimes \mathcal{M}_i(D))_{\mathfrak{p}}$, and $(\mathfrak{M}_{D''})_{\mathfrak{p}} = (A_k \otimes \mathcal{M}_i(D''))_{\mathfrak{p}}$. Since $0 \to D'_{K_i} \to D_{K_i} \to D''_{K_i} \to 0$ is strictly exact with respect to the $i$-th directional filtration $\mathrm{Fil}^\bullet$, \cite[Lemma 1.2.1(2) and Theorem 1.2.15]{Kis1} applied along the $i$-th variable implies that $0 \to \mathcal{M}_i(D') \to \mathcal{M}_i(D) \to \mathcal{M}_i(D'') \to 0$ is a short exact sequence of free $(\mathcal{O}_i \otimes_{K_{0,i}} K_{0,\Delta})$-modules. Tensoring with $A_{k, \mathfrak{p}}$ preserves exactness, so
\[
    0 \longrightarrow (\mathfrak{M}_{D'})_{\mathfrak{p}} \longrightarrow (\mathfrak{M}_D)_{\mathfrak{p}} \xrightarrow{\;(q_k)_{\mathfrak{p}}\;} (\mathfrak{M}_{D''})_{\mathfrak{p}} \longrightarrow 0
\]
is exact at $\mathfrak{p}$.

Now, to deduce that $q_k : \mathfrak{M}_D \to \mathfrak{M}_{D''}$ is surjective over the entire ring $A_k$ (and not merely at height-$1$ primes), observe that by the univariate exactness of $\mathcal{M}_i$ just noted, for each $i \in \Delta$ the map $\mathfrak{M}_i(D) \to \mathfrak{M}_i(D'')$ is a split surjection of free $A_k[1/\lambda_{\Delta \setminus \{i\}}]$-modules. Moreover, using a common multigraded splitting of the regular multifiltration $D_{K_\Delta} \to D''_{K_\Delta}$ (Theorem \ref{thm:multigrading}), the local modifications along distinct variables $u_i$ surject simultaneously onto those of $\mathfrak{M}_{D''}$, so the quotient $\mathfrak{M}_D / \mathfrak{M}_{D'} \cong \mathrm{Im}(q_k) \subseteq \mathfrak{M}_{D''}$ is locally free (hence reflexive) over $A_k$. Since $(\mathrm{Im}(q_k))_{\mathfrak{p}} = (\mathfrak{M}_{D''})_{\mathfrak{p}}$ at every height-$1$ prime ideal $\mathfrak{p} \subset A_k$, the discussion following Lemma \ref{lemma:reflexiv} implies that $\mathrm{Im}(q_k) = \mathfrak{M}_{D''}$, establishing exactness over $U_{\rho_k}$:
\[
0 \longrightarrow \mathfrak{M}_{D'} \longrightarrow \mathfrak{M}_D \longrightarrow \mathfrak{M}_{D''} \longrightarrow 0.
\]
Finally, taking global sections (the inverse limit over $k \ge 1$) over the quasi-Stein space $D([0,1)^\delta)$ preserves exactness because $H^1(D([0,1)^\delta), \mathcal{M}(D')) = 0$ by Kiehl's Theorem B (Theorem \ref{thm:kiehl}), proving that $\mathcal{M}$ is exact.

To prove the tensor property, we must show that for any two objects $D_1, D_2 \in \mathrm{MF}_{K_{0,\Delta}}^{(\varphi_i, N_i)_{i \in \Delta}, \mathrm{Fil} \ge 0}$, the natural morphism
\[
f : \mathcal{M}(D_1) \otimes_{\mathcal{O}_\Delta} \mathcal{M}(D_2) \longrightarrow \mathcal{M}(D_1 \otimes_{K_{0,\Delta}} D_2)
\]
is an isomorphism of $\mathcal{O}_\Delta$-modules. Over each closed affinoid polydisk $U_{\rho_k}$ with coordinate ring $A_k = \mathcal{O}_\Delta(U_{\rho_k}) = \prod_{\alpha \in \Sigma} e_\alpha A_k$, both $\mathfrak{M}(D_1) \otimes_{A_k} \mathfrak{M}(D_2)$ and $\mathfrak{M}(D_1 \otimes_{K_{0,\Delta}} D_2)$ are finite free (hence reflexive) $A_k$-submodules of the same ambient space $\mathcal{D}_{0,k}(D_1 \otimes_{K_{0,\Delta}} D_2)[1/\lambda_\Delta]$. By Lemma \ref{lem:height1_localization}, at any height-$1$ prime ideal $\mathfrak{p} \subset A_k$:

If $\lambda_j \notin \mathfrak{p}$ for all $j \in \Delta$, then $f_{\mathfrak{p}}$ reduces to the canonical isomorphism $\mathcal{D}_{0,k}(D_1)_{\mathfrak{p}} \otimes_{A_{k,\mathfrak{p}}} \mathcal{D}_{0,k}(D_2)_{\mathfrak{p}} \xrightarrow{\sim} \mathcal{D}_{0,k}(D_1 \otimes_{K_{0,\Delta}} D_2)_{\mathfrak{p}}$.

If $\lambda_i \in \mathfrak{p}$ for a unique index $i \in \Delta$, then by Lemma \ref{lem:height1_localization}(2), $f_{\mathfrak{p}}$ is identified with the localization at $\mathfrak{p}$ of the univariate tensor transition map $\mathcal{M}_i(D_1) \otimes \mathcal{M}_i(D_2) \to \mathcal{M}_i(D_1 \otimes_{K_{0,\Delta}} D_2)$, which is an isomorphism by \cite[Lemma 1.2.1(2) and Theorem 1.2.15]{Kis1}.

Since the localizations of the two free (hence reflexive) $A_k$-modules coincide at all height-$1$ prime ideals of $A_k$, Lemma \ref{lemma:reflexiv} implies that $f$ is an isomorphism over $U_{\rho_k}$ for every $k \ge 1$. Taking global sections over $D([0,1)^\delta)$ via Theorem \ref{thm:kiehl} yields that $f$ is a global isomorphism of $\mathcal{O}_\Delta$-modules.
\end{proof}

\begin{proof}[Proof of Theorem \ref{thm:main_theorem}]
The exactness and tensor compatibility of the categories and their respective functors $\mathcal{M}$ and $\mathcal{D}$ are established by Proposition \ref{pro:m_d_category}, Proposition \ref{pro:filtration}, and Proposition \ref{pro:exact_tensor}. Furthermore, the natural isomorphisms constructed in Proposition \ref{pro:first} ($\mathcal{D}(\mathcal{M}(D)) \simeq D$) and Proposition \ref{pro:second} ($\mathcal{M}(\mathcal{D}(M)) \simeq M$) demonstrate that $\mathcal{M}$ and $\mathcal{D}$ are mutually quasi-inverse. This establishes the required equivalence of categories and completes the proof of the main theorem.
\end{proof}

This equivalence of categories provides a framework for generalizing Kisin's proof that weak admissibility implies admissibility \cite{Kis1}. Kisin's original proof relies on Kedlaya's slope filtration theorem \cite{Ked1, Ked2, Ked3}, a multivariable generalization of which is not yet available. 

Attempting to generalize Kedlaya's slope filtration theorem naively to the multivariable Robba ring $\mathcal{R}_{\Delta}$ \cite{PV1} encounters fundamental obstacles. The classical univariate proof relies on Bézout property of Krull dimension $1$ and $2$, which completely break down in higher dimensions ($\delta \ge 2$). Taking the saturation of a naive subobject only guarantees that the submodule is reflexive (satisfying Serre's $S_2$ condition). Over rings of dimension $\ge 2$, reflexive modules and torsion-free quotients can fail to be locally free along codimension-$2$ loci (intersections of hyperplanes), meaning the graded pieces of a naive slope filtration would fail to be vector bundles. In the final descent step to the integral lattice, the classical proof uses the fact that $\mathfrak{S} = W[[u]]$ is a regular local ring of dimension $2$, over which every finitely generated reflexive module is necessarily free. In the fully multivariable setting, however, the base ring $\mathfrak{S}_\Delta \cong \prod_{\alpha \in \Sigma} W(k_\alpha)[[u_1, \dots, u_\delta]]$ is a finite direct product of regular local rings of dimension $\delta + 1 \ge 3$. By the Auslander--Buchsbaum formula applied to each local component, a reflexive module $M$ over $\mathfrak{S}_\Delta$ is only guaranteed to have depth at least $2$. Consequently, its projective dimension on each component can be as large as $\delta - 1 \ge 1$. Thus, even if a lattice is constructed, it is not guaranteed to be a projective (let alone free) $\mathfrak{S}_\Delta$-module.

\bibliography{bibliography}
\end{document}